\documentclass[12pt]{amsart}
\usepackage{amssymb}
\theoremstyle{plain}
  \newtheorem{thm}{Theorem}[section]
  \newtheorem{lem}[thm]{Lemma}
  
  \newtheorem{cor}[thm]{Corollary}
  \newtheorem{prop}[thm]{Proposition}
  \newtheorem{clm}[thm]{Claim}

\theoremstyle{definition}
  \newtheorem{defn}[thm]{Definition}
  \newtheorem{conj}[thm]{Conjecture}
  \newtheorem{ex}[thm]{Example}
  
  \newtheorem{prob}[thm]{Problem}
  
  \newtheorem{nota}[thm]{Notation}

\theoremstyle{remark}
  \newtheorem{rem}[thm]{Remark}

\newenvironment{proclaim}[1]
  {\par\medskip\noindent #1 \ }
  {\par\medskip}

\numberwithin{equation}{section}

\DeclareMathOperator{\supp}{supp}

\DeclareMathOperator{\Dom}{\mathcal{D}}

\DeclareMathOperator{\CAT}{CAT}
\DeclareMathOperator{\ord}{ord}
\DeclareMathOperator{\dis}{dis}

\newcommand{\calE}{\mathcal{E}}

\newcommand{\field}[1]{\mathbb{#1}}

\newcommand{\C}{\field{C}}
\newcommand{\R}{\field{R}}

\newcommand{\N}{\field{N}}

\renewcommand{\liminf}{\mathop{\varliminf}}
\renewcommand{\limsup}{\mathop{\varlimsup}}

\begin{document}
\title[Variational convergence over metric spaces]
{Variational convergence\\ over metric spaces\\
  \footnote{\today}}

\author{Kazuhiro Kuwae}
\address{Department of Mathematics, Faculty of Education,
  Kumamoto University, Kumamoto, 860-8555, JAPAN}
\email{kuwae@gpo.kumamoto-u.ac.jp}

\author{Takashi Shioya}
\address{Mathematical Institute, Tohoku University, Sendai 980-8578,
  JAPAN}
\email{shioya@math.tohoku.ac.jp}





\thanks{The first author is partially supported by a Grant-in-Aid
  for Scientific Research No.~16540201 from the Ministry of
  Education, Science, Sports and Culture, Japan}

\thanks{The second author is partially supported by a Grant-in-Aid
  for Scientific Research No.~14540056 from the Ministry of
  Education, Science, Sports and Culture, Japan}

\begin{abstract}
  We introduce a natural definition of $L^p$-convergence of maps, $p
  \ge 1$, in the case where the domain is a convergent sequence of
  measured metric space with respect to the measured Gromov-Hausdorff
  topology and the target is a Gromov-Hausdorff convergent sequence.
  With the $L^p$-convergence, we establish a theory of variational
  convergences.  We prove that the Poincar\'e inequality with some
  additional condition implies the asymptotic compactness.  The
  asymptotic compactness is equivalent to the Gromov-Hausdorff
  compactness of the energy-sublevel sets.  Supposing that the targets
  are $\CAT(0)$-spaces, we study convergence of resolvents.  As
  applications, we investigate the approximating energy functional
  over a measured metric space and convergence of energy functionals
  with a lower bound of Ricci curvature.
\end{abstract}

\maketitle

\tableofcontents

\section{Introduction}

Let $M_i \to M$ and $Y_i \to Y$ ($i = 1,2,3,\dots$) be two pointed
Gromov-Hausdorff convergent sequences of proper metric spaces, where
`\emph{proper}' means that any bounded subset is relatively compact,
and let us give measures on $M_i$ which converge to a measure on $M$.
We are interested in the convergence and asymptotic behavior of maps
$u_i : M_i \to Y_i$ and also energy functionals $E_i$ defined on a
family of maps from $M_i \to Y_i$.  Note that there are several
attempts to define natural energy functionals on the mapping space
from $M$ to $Y$ by the measured metric structure of $M$ and the metric
structure of $Y$ (see, for example, \cite{HK:SmetP, St:heat,
  Ch:metmeas} for the case of $Y = \R$ and \cite{Jt:equilib,
  KvS:sobharm, KwSy:sobmet, Ota:Cheeger} for the general $Y$).  We
introduce a natural definition of the $L^p$-convergence of $u_i : M_i
\to Y_i$ to $u : M \to Y$, $p \ge 1$, and establish a general theory
for energy functionals $E_i$ by extending the theory of variational
convergences, mainly studied by Mosco \cite{Ms:compmedia}.  Mosco
introduced the asymptotic compactness of energy functionals $\{E_i\}$,
which is a generalization of the Rellich compactness.  The asymptotic
compactness is useful to obtain the convergence of energy minimizers,
i.e., harmonic maps, and also to investigate spectral properties in
the linear case.  Under a uniform bound of Poincar\'e constants and
some property of the metric on $M$, we prove the asymptotic
compactness of $\{E_i\}$.  We focus on a $\Gamma$-convergence with the
asymptotic compactness, say \emph{compact convergence}.  If $\{E_i\}$
is asymptotically compact, it has a compact convergent subsequence.
We prove that the compact convergence is equivalent to the
Gromov-Hausdorff convergence of the energy-sublevel sets, which is a
geometric interpretation of compact convergence.  We also prove that
the compact convergence of $E_i$ is equivalent to a convergence of
associated resolvents, provided that $Y_i$ are all $\CAT(0)$-spaces
and $E_i$ are lower semi-continuous convex functionals.  Such the
resolvent was defined by Jost \cite{Jt:convex} using the Moreau-Yosida
approximation.  As applications of our theory, we study the
approximating energy functional and its spectral property.  We also
obtain the compactness of energy functionals if $M_i$ are Riemannian
manifolds with a lower bound of Ricci curvature.

We mention further details.  Let $M_i \to M$ and $Y_i \to Y$ be as
above.  For $p \ge 1$ and two measurable maps $u, v : M \to Y$, the
\emph{$L^p$-distance $d_{L^p}(u,v)$ between $u$ and $v$} is defined by
\[
d_{L^p}(u,v) := \left( \int_M d_Y(u(x),v(x))^p \; dx
\right)^{\frac{1}{p}}.
\]
Consider the $L^p$-metric space of measurable maps from $M$ to $Y$
(resp. from $M_i$ to $Y_i$) and denote it by $X$ (resp.~$X_i$).  To
define the $L^p$-convergence of maps $u_i \in X_i$ to a map $u \in X$,
we take Gromov-Hausdorff approximations $\varphi_i : M_i \to M$, which
are (not necessarily continuous) Borel maps almost preserving the
distances.  Assume first that $Y_i = Y$ for all $i$.  It is natural to
consider that $\varphi_i \circ u$ converges to $u$, so that if a
sequence $u_i : M_i \to Y (= Y_i)$ satisfies $d_{L^p}(\varphi_i \circ
u,u_i) \to 0$ as $i \to \infty$, then we say that \emph{$u_i$
  $L^p$-converges to $u$}.  However, this definition is natural only
if $u$ is continuous, because the convergence of the measures on $M_i$
is only weak (see Remark \ref{rem:Lptop} for the detailed
explanation).  If $Y_i = Y$ is (a subset of) a Banach space, then we
can take a continuous $L^p$-approximation $\tilde u_\epsilon$ of a
measurable $u : M \to Y$, $\epsilon > 0$, such that $d_{L^p}(\tilde
u_\epsilon,u) < \epsilon$.  If
\begin{equation}
  \label{eq:Lpdef}
  \lim_{\epsilon \to 0} \limsup_{i \to \infty} d_{L^p}(\varphi_i \circ
  \tilde u_\epsilon,u_i) = 0,
\end{equation}
then $u_i$ $L^p$-converges to $u$.  For general $Y_i$, $Y$, we embed all
$Y_i$ and $Y$ into a Banach space and employ the same definition as
above.  This induces a topology on the disjoint union $\bigsqcup_i X_i
\sqcup X$, say the \emph{$L^p$-topology}.
We prove that the $L^p$-topology is independent of
the Gromov-Hausdorff approximation $\{\varphi_i\}$ and
also of the embedding of $Y_i, Y$ into a Banach space.
See Section \ref{ssec:Lp} for the $L^p$-topology.

The $L^p$-topology on $\bigsqcup_i X_i \sqcup X$ has some nice
properties involving the $L^p$-metric structure of $X_i$ and $X$, such
as, if $X_i \ni u_i,v_i \to u,v \in X$ respectively in $L^p$, then
$d_{L^p}(u_i,v_i) \to d_{L^p}(u,v)$.  By
their properties we present a set of axioms for a topology on
$\bigsqcup_i X_i \sqcup X$ for general metric spaces $(X_i,d_{X_i})$
and $(X,d_X)$.  We call such a topology satisfying the axioms the
\emph{asymptotic relation between $\{X_i\}$ and $X$} (see Definition
\ref{defn:asymprel}).  Since $X_i$ and $X$ are typically improper, the
asymptotic relation can be thought as a non-uniform variant of
Gromov-Hausdorff convergence $X_i \to X$.

After establishing the foundation of asymptotic relation, we give, in
Section \ref{ssec:def}, a formulation of variational convergences of
general functions $E_i : X_i \to [\,0,+\infty\,]$ to $E : X \to
[\,0,+\infty\,]$.  $\{E_i\}$ is said to be \emph{asymptotically
  compact} if for any bounded sequence $u_i \in X_i$ with $\sup_i
E_i(u_i) < +\infty$, it has a convergent subsequence, where
`\emph{bounded}' means that $d_{X_i}(u_i,v_i)$ is bounded for some
convergent sequence $v_i \in X_i$.  We say that \emph{$E_i$
  $\Gamma$-converges to $E$} if the following ($\Gamma$1) and
($\Gamma$2) are both satisfied.
\begin{enumerate}
\item[($\Gamma$1)] For any $u \in X$ there exists a sequence $u_i \in
  X_i$ such that $u_i \to u$ and $E_i(u_i) \to E(u)$.
\item[($\Gamma$2)] We have $\liminf_i E_i(u_i) \ge E(u)$ for any
  convergent sequence $X_i \ni u_i \to u \in X$.
\end{enumerate}
We say that \emph{$E_i$ compactly converges to $E$} if $E_i$
$\Gamma$-converges to $E$ and if $\{E_i\}$ is asymptotically compact.

In Section \ref{sec:varconvmet}, we study the asymptotic compactness
and the compact convergence of $\{E_i\}$.  It is important to
investigate under what condition the asymptotic compactness is
obtained.  We introduce a concept of the \emph{local covering order}
of a locally compact metric space (see Definition \ref{defn:loccov}),
which is a quantity measuring the local size of the metric space.
Assume that $M_i$ and $M$ are all compact.  Under a bound of local
covering order of $M$ and a uniform bound of Poincar\'e constants for
$E_i$ on $M_i$, we prove the asymptotic compactness of $\{E_i\}$ (see
Theorem \ref{thm:Pasympcpt}).  Since we do not need the doubling
condition, our theorem can be applied to infinite-dimensional spaces.
There is a mistake in the linear version in \cite{KwSy:specstr}.  We
correct it by introducing the local covering order.

Motivated by Gromov's study of spectral concentration, Section
3$\tfrac12$.57 of \cite{Gr:greenbook}, we prove that the compact
convergence $E_i \to E$ is equivalent to that for any $c \in \R$ there
exists a sequence $c_i \searrow c \in \R$ such that the sublevel set
$(\{\,u \in X_i \mid E_i(u) \le c_i\,\},o_i)$ converges to $(\{\,u \in
X \mid E(u) \le c\,\},o)$ with respect to the pointed Gromov-Hausdorff
topology, where $o_i \in X_i$ is a sequence converging to a point $o
\in X$ (see Theorem \ref{thm:cpt}).

In Section \ref{sec:CAT0}, we study variational convergences over
$\CAT(0)$-spaces, where a $\CAT(0)$-space is a globally nonpositively
curved metric space.  Typical examples of $\CAT(0)$-spaces are
Hadamard manifolds and trees.  If the target space is $\CAT(0)$ then
the $L^2$-mapping space is also $\CAT(0)$, and an energy functional
defined in a suitable way becomes convex and lower semi-continuous.
Thus, it is reasonable to assume that $X_i$ and $X$ are all
$\CAT(0)$-spaces, and $E_i$ and $E$ are convex lower semi-continuous
functions with $E_i, E \not\equiv +\infty$.  For any $\lambda \ge 0$
and $u \in X$, there exists a unique minimizer, say $J^E_\lambda(u)
\in X$, of $v \mapsto \lambda E(v) + d_X(u,v)^2$.  This defines a map
$J^E_\lambda : X \to X$, called the \emph{resolvent of $E$} (see
\cite{Jt:convex}).  Note that if $X$ is a Hilbert space and if $E$ is
a closed densely defined symmetric quadratic form on $X$, then we have
$J^E_\lambda = (I + \lambda A)^{-1}$, where $A$ is the infinitesimal
generator associated with $E$.  The one-parameter family
$[\,0,+\infty\,) \ni \lambda \mapsto J^E_\lambda(u)$ gives a
deformation of a given map $u \in X$ to a minimizer of $E$ (or a
harmonic map), $\lim_{\lambda \to +\infty} J^E_\lambda(u)$ (if any).
Jost \cite{Jt:nlinDir} studied convergence of resolvents and
Moreau-Yosida approximations.  Although his study is only on a fixed
$\CAT(0)$-space, we extend it for a sequence of $\CAT(0)$-spaces with
an asymptotic relation.  We investigate the relation between
asymptotic compactness and resolvents (see Proposition
\ref{prop:EJasympcpt}).  This is new even on a fixed $\CAT(0)$-space.
As the main theorem in this section (see Theorem \ref{thm:cptconv}),
we prove that $E_i$ compactly converges to $E + c$ for some constant
$c \in \R$ iff the following (1) and (2) are both satisfied.
\begin{enumerate}
\item For any $\lambda > 0$ and for any bounded sequence $u_i \in
  X_i$, $J^{E_i}_\lambda(u_i)$ has a convergent subsequence.
\item $J^{E_i}_\lambda(u_i)$ converges to $J^E_\lambda(u)$ for any
  $\lambda > 0$ and for any convergent sequence $X_i \ni u_i \to u \in
  X$.
\end{enumerate}
In the case where $E_i$ and $E$ are symmetric quadratic forms on
Hilbert spaces, this theorem is obtained in our previous paper
\cite{KwSy:specstr} (see also Section \ref{ssec:Hil} of this paper)
and is applied to the study of the spectral properties of the
Laplacian of convergent Riemannian manifolds.  This also has some
applications to a homogenization problem, \cite{V:sound}, and
convergence of Dirichlet forms, \cite{Ks:convLapII, Ks:graph,
  OST:conv}, including finite-dimensional approximation problems,
\cite{Kol:convDir, Kol:convinf}.  The results of such studies could
possibly extend to the nonlinear case.

In Section \ref{sec:appl}, we give some applications.  Define the
\emph{$\rho$-approximating energy $E^\rho(u)$, $\rho > 0$, of a
  measurable map $u : M \to Y$} by
\[
E^\rho(u) := \frac{1}{2} \int_M \frac{1}{|B(x,\rho)|}
\int_{B(x,\rho) \setminus \{x\}}
\frac{d_Y(u(x),u(y))^p}{\rho^p}\;dydx \in [\,0,+\infty\,],
\]
where $|B(x,\rho)|$ is the measure of the open metric $\rho$-ball
$B(x,\rho)$ centered at $x$.  If $M$ and $Y$ are both complete
Riemannian manifolds, then as $\rho \to 0+$, $E^\rho$
$\Gamma$-converges to the usual energy functional upto a constant
multiple over the $L^p$-mapping space.  If $M$ is a Riemannian
manifold and if $Y$ is a metric space, then $E^\rho$
$\Gamma$-converges to some functional $E$, which is often called the
\emph{Korevaar-Schoen energy functional} (see \cite{KvS:sobharm}).  We
prove that $\{E^\rho\}$ is asymptotically compact and so the
convergence $E^\rho \to E$ becomes compact, provided $M$ is a compact
Riemannian manifold and $Y$ is a proper metric space.  More generally,
this is true in the case where $M$ is a compact measured metric space
with some property (called the measure contraction property) stated in
\cite{St:heat, KwSy:sobmet}.  Note that each $E^\rho$ does not have
the Rellich compactness property.  This brings an example of an
asymptotically compact sequence $\{E_i\}$ with noncompact $E_i$.  In
the real-valued case (i.e., $Y = \R$), the functionals $E^\rho$ and
$E$ are symmetric quadratic forms and we denote their infinitesimal
generators by $A^\rho$ and $A$ respectively.  The operator $A$ is the
Laplacian if $M$ is a Riemannian manifold.  Each $A^\rho$ has nonempty
essential spectrum, whose bottom is divergent to $+\infty$ as $\rho
\to 0+$.  The $k^{th}$ eigenvalue of $A^\rho$ converges to that of $A$
as $\rho \to 0+$ for any fixed $k \in \N$.

We next study the convergence of energy functionals under a lower
bound of Ricci curvature.  For a constant $n \ge 2$, let $M_i$ be a
sequence of $n$-dimensional closed Riemannian manifolds of Ricci
curvature $\ge -(n-1)$ which converges to a compact measured metric
space $M$ with respect to the measured Gromov-Hausdorff topology.  Let
$Y_i \to Y$ be a pointed Gromov-Hausdorff convergent sequence of
proper pointed metric spaces.  Consider the Korevaar-Schoen type
energy, denote $E_i$, on the $L^2$-mapping space from $M_i$ to $Y_i$
(see \cite{KvS:sobharm}).  We prove that $\{E_i\}$ is asymptotically
compact and has a compactly convergent subsequence (see Theorem
\ref{thm:Ric}).  If $Y_i = Y$ is a fixed complete Riemannian manifold,
then Kasue's result \cite{Ks:convLapII} together with
Cheeger-Colding's one \cite{CC:strRicIII} shows that $E_i$ converges
in some sense to a naturally defined energy functional $E$ on the
$L^2$-mapping space from $M$ to $Y$.  Combining this with our result
yields that $E_i$ compactly converges to $E$.

\section{Preliminaries for Gromov-Hausdorff convergence}
\label{sec:prelim}

The \emph{Hausdorff distance} $d_H^Z(X,Y)$ between two subsets $X$ and
$Y$ of a metric space $Z$ is defined to be the infimum of $r > 0$
such that $X \subset B(Y,r)$ and $Y \subset B(X,r)$, where $B(X,r)$
denotes the open $r$-neighborhood of $X$.  Let $X$ and $Y$ be two
compact metric spaces.  The \emph{Gromov-Hausdorff distance}
$d_{GH}(X,Y)$ between $X$ and $Y$ is the infimum of $d_H^Z(X,Y)$,
where $Z$ is any metric space into which $X$ and $Y$ are isometrically
embedded.  The \emph{distortion} $\dis\varphi$ of a (not necessarily
continuous) map $\varphi : X \to Y$ is defined by
\[
\dis \varphi := \sup_{x,y \in X} | d(\varphi(x),\varphi(y)) - d(x,y) |,
\]
where $d$ denotes the distance function.  If a map $\varphi : X \to Y$
satisfies $\dis\varphi < \epsilon$ and $B(\varphi(X),\epsilon) = Y$,
then it is called an \emph{$\epsilon$-approximation}.  It is known
that if $d_{GH}(X,Y) < \epsilon$ then there exists a
$2\epsilon$-approximation from $X$ to $Y$.  Conversely, if there
exists an $\epsilon$-approximation from $X$ to $Y$ then $d_{GH}(X,Y) <
2\epsilon$.

Let $\{i\}$ be a directed set and $\{(X_i,o_i)\}$ a net of pointed
proper metric spaces, where `\emph{proper}' means that any closed and
bounded subset is compact.  We say that \emph{$(X_i,o_i)$ converges to
  a pointed proper metric space $(X,o)$ with respect to the pointed
  Gromov-Hausdorff convergence} if for any $r > 0$ there exist two
nets of positive numbers $r_i \searrow r$, $\epsilon_i \to 0+$, and
$\epsilon_i$-approximations $\varphi_i : B(o_i,r_i) \to B(o,r)$ with
$\varphi_i(o_i) = o$.  This is equivalent to the existence of
$\epsilon_i$-approximations $\psi_i : B(o_i,r_i') \to B(o,r_i)$ such
that $\psi_i(o_i) = o$, $\epsilon_i \to 0+$ and $r_i' > r_i \to
+\infty$.  We call such a $\{\psi_i\}$ a \emph{pointed
  Gromov-Hausdorff approximation}.  This convergence induces a
topology, called the \emph{pointed Gromov-Hausdorff topology}, on the
set of pointed proper metric spaces.

Note that a notion of convergence of a sequence of countable elements
is not enough to define a topology and we need convergence of a net
with indexed by a directed set for it.

The following definition seems not to be in a common knowledge.

\begin{defn}[Compact Hausdorff convergence]
  Let $A_i, A$ be closed subsets of a metric space $Z$.  We say that
  \emph{$A_i$ converges to $A$ in the compact Hausdorff convergence}
  if for some point $o \in Z$ and for any $r > 0$ there
  exists a net $r_i \searrow r$ such that
  \[
  d_H^Z(A_i \cap B(o,r_i),A \cap B(o,r)) \to 0.
  \]
\end{defn}

\begin{rem}
  If the restrictions of the metric on $A_i$ and $A$ are all proper,
  then a compact Hausdorff convergence $A_i \to A$ implies the pointed
  Gromov-Hausdorff convergence $(A_i,o_i) \to (A,o)$ for any net $o_i
  \in X$ converging to a point $o \in X$.
\end{rem}

The following proposition seems to be well-known.

\begin{prop}[cf.~Proof of 3.5(b) of \cite{Gr:greenbook}]
\label{prop:GHunion}
  Let $(X_i,o_i)$ be a net of pointed proper metric spaces converging
  to a pointed proper metric space $(X,o)$ in the pointed Gromov-Hausdorff
  topology.  Then, there exists a proper metric $d_\mathcal{X}$ on the
  disjoint union $\mathcal{X} := \bigsqcup_i X_i \sqcup X$ such that
  \begin{enumerate}
  \item the restrictions of $d_\mathcal{X}$ on $X_i$ and $X$ coincide with the
    original metrics $d_{X_i}$ and $d_X$ respectively;
  \item $X_i$ converges to $X$ in the compact Hausdorff convergence in
    $(\mathcal{X},d_\mathcal{X})$;
  \item $d_{\mathcal{X}}(o_i,o) \to 0$.
  \end{enumerate}
\end{prop}

\begin{proof}
  Since the proof is standard, we here give an outline.  Take a dense
  countable subset $\{p_n\}_{n \in \N} \subset X$ with $p_1 = 0$.  By
  the convergence $(X_i,o_i) \to (X,o)$, there exist a net $N(i)
  \nearrow \infty$ of natural numbers and points $p_{n,i} \in X_i$ for
  $n=1,2,\dots,N(i)$ with $p_{1,i} = o_i$ such that
  \[
  |\;d_{X_i}(p_{m,i},p_{n,i}) - d_X(p_m,p_n)\;| < 1/N(i)
  \]
  for any $m,n = 1,2,\dots,N(i)$.  Set $d_\mathcal{X} := d_{X_i}$ on $X_i
  \times X_i$ and $d_\mathcal{X} := d_X$ on $X \times X$.  We define the
  distance $d_\mathcal{X}(x,y)$ between any two points $x \in X_i$ and $y \in
  X$ by
  \[
  d_\mathcal{X}(x,y) := \inf_{n=1,2,\dots,N(i)}
  d_{X_i}(x,p_{n,i}) + d_X(y,p_n) + 1/N(i).
  \]
  and the distance $d_\mathcal{X}(x,y)$ between $x \in X_i$ and $y \in X_j$
  for $i \neq j$ by
  \[
  d_\mathcal{X}(x,y) := \inf_{z \in X} d_\mathcal{X}(x,z) + d_\mathcal{X}(y,z).
  \]
  Then, $d_\mathcal{X}$ is a unique minimal distance function on
  $\mathcal{X}$ satisfying (1) and the condition that
  $d_\mathcal{X}(p_{n,i},p_n) = 1/N(i)$ for all $n = 1,2,\dots,N(i)$.
  In particular we have (3).  It is easy to verify (2).  This
  completes the proof.
\end{proof}

We assume that all measure spaces are locally compact Polish spaces
with positive Radon measures of full support, where a \emph{Polish
  space} is, by definition, homeomorphic to a complete separable
metric space.

\begin{defn}[Measure approximation, \cite{KwSy:specstr}]
  \label{defn:mapprox}
  Let $M_i$ and $M$ be measure spaces.  A net $\{\varphi_i : M_i
  \supset \Dom(\varphi_i) \to M\}$ of maps is called a \emph{measure
    approximation} if the following (M1) and (M2) are satisfied.
  \begin{itemize}
  \item[(M1)] Each $\varphi_i$ is a measurable map from a Borel subset
    $\Dom(\varphi_i)$ of $M_i$ to $M$.
  \item[(M2)] The push-forward by $\varphi_i$ of the measure on $M_i$
    \emph{vaguely} (or \emph{weakly-star}) converges to the measure on
    $M$, i.e., for any $f \in C_0(M)$,
  \[
  \lim_i \int_{\Dom(\varphi_i)} f \circ \varphi_i(x) \; dx
  = \int_M f(x) \; dx,
  \]
  where $C_0(M)$ is the set of continuous functions on $M$ with
  compact support.
  \end{itemize}
\end{defn}

Note that it is nonsense to define a topology on the set of measure
spaces by the existence of a measure approximation.  This is because,
if the total measure of $M_i$ converges to a finite number $a$, then
the map from $M_i$ to the set of a single point with mass $a$ forms a
measure approximation.  However, the following definition makes
sense.

\begin{defn}[Measured Gromov-Hausdorff topology, \cite{Fk:laplace,
  Gr:greenbook}]
  Let $M_i$ and $M$ be compact measured metric spaces.  We say that
  \emph{$M_i$ converges to $M$ in the sense of the measured
    Gromov-Hausdorff convergence} if there exists a measure
  approximation $\{\varphi_i : M_i \to M\}$ such that each $\varphi_i$
  is an $\epsilon_i$-approximation for some $\epsilon_i \to 0+$.
  
  The pointed version is also defined.  Let $M_i$ and $M$ be pointed
  proper measured metric spaces.  We say that \emph{$M_i$ converges to
    $M$ in the sense of the pointed measured Gromov-Hausdorff
    convergence} if there exists a pointed Gromov-Hausdorff
  approximation that is a measure approximation, which we call a
  \emph{measured pointed Gromov-Hausdorff approximation}.
  
  The (pointed) measured Gromov-Hausdorff convergence induces a
  topology, called the (\emph{pointed}) \emph{measured
    Gromov-Hausdorff topology}, on the set of compact (resp.~pointed
  proper) measured metric spaces, which is stronger than the
  (pointed) Gromov-Hausdorff topology.
\end{defn}

\section{Asymptotic relation} \label{sec:asymprel}

In Section \ref{ssec:found}, we present a set of axioms of asymptotic
relation and prove some lemmas needed in the later sections.  Then,
in Section \ref{ssec:Lp}, we prove that $L^p$-mapping spaces $X_i$ and
$X$ as in Introduction satisfy the axioms.  The axioms are more
flexible than our previous ones in the linear case,
\cite{KwSy:specstr}.

\subsection{Foundation} \label{ssec:found}

Throughout this paper, we denote by $i$ any element of a given directed set
$\{i\}$.  Let $\{X_i\}$ be a net of metric spaces and $X$ a metric space.
Define
\[
\mathcal{X} := \left( \bigsqcup_i X_i \right) \sqcup X
\qquad (\text{disjoint union}).
\]
We sometimes consider the following additional condition:
\begin{itemize}
\item[(L)] Given a field $K = \R$ or $\C$, $X_i$ and $X$ are all
  topological linear spaces over $K$ whose topologies are compatible
  with their metric structure.
\end{itemize}

\begin{defn}[Asymptotic relation] \label{defn:asymprel}
  We call a topology on $\mathcal{X}$ satisfying the following (A1)--(A4) an
  \emph{asymptotic relation between $\{X_i\}$ and $X$}.
  \begin{itemize}
  \item[(A1)] $X_i$ and $X$ are all closed in $\mathcal{X}$ and the
    restricted topology of $\mathcal{X}$ on each of $X_i$ and $X$ coincides
    with its original topology.
  \item[(A2)] For any $x \in X$ there exists a net $x_i \in X_i$
    converging to $x$ in $\mathcal{X}$.
  \item[(A3)] If $X_i \ni x_i \to x \in X$ and $X_i \ni y_i \to y \in
    X$ in $\mathcal{X}$, then we have $d_{X_i}(x_i,y_i) \to d_X(x,y)$.
  \item[(A4)] If $X_i \ni x_i \to x \in X$ in $\mathcal{X}$ and if $y_i \in
    X_i$ is a net with $d_{X_i}(x_i,y_i) \to 0$, then $y_i \to x$ in
    $\mathcal{X}$.
  \end{itemize}
  Supposing (L), we say that an asymptotic relation between $\{X_i\}$
  and $X$ is \emph{linear} if the following (AL) is satisfied.
  \begin{itemize}
  \item[(AL)] If $X_i \ni x_i \to x \in X$ and $X_i \ni y_i \to y \in
    X$ in $\mathcal{X}$, then $ax_i+by_i \to ax+by$ in $\mathcal{X}$ for any
    scalers $a, b \in K$.
  \end{itemize}
\end{defn}

By the definition, an asymptotic relation on $\mathcal{X}$ is a
Hausdorff topology.

\begin{rem}
  Notice that a Gromov-Hausdorff convergence $X_i \to X$ induces an
  asymptotic relation (see Proposition 3.5(b) of \cite{Gr:greenbook}
  and also Proposition \ref{prop:GHunion} of this paper).  Thus, the
  existence of an asymptotic relation between $\{X_i\}$ and $X$ can be
  thought as a generalization of Gromov-Hausdorff convergence $X_i \to
  X$ in a sense.  However, it seems not to be suitable to call an
  asymptotic relation a convergence.  This is because the existence of
  an asymptotic relation between $\{X_i\}$ and $X$ induces one
  between $\{X_i\}$ and any $Y \subset X$ as the restriction.
\end{rem}

\begin{defn}[Metric approximation]
  A net $\{f_i : X \supset \Dom(f_i) \to X_i\}$ of (not necessarily
  continuous) maps is called a \emph{metric approximation for
    $\{X_i\}$ and $X$} if the following (B1) and (B2) are satisfied.
  \begin{itemize}
  \item[(B1)] $\Dom(f_i)$ is monotone nondecreasing in $i$ and
    $\bigcup_i \Dom(f_i)$ is dense in $X$.
  \item[(B2)] For any $x,y \in \bigcup_i \Dom(f_i)$ we have
    $d_{X_i}(f_i(x),f_i(y)) \to d_X(x,y)$.
  \end{itemize}
  Under (L), a metric approximation $\{f_i : X \supset \Dom(f_i) \to
  X_i\}$ is said to be \emph{linear} if the following (BL) is satisfied.
  \begin{itemize}
  \item[(BL)] Each $f_i$ is a linear map from a linear subspace
    $\Dom(f_i) \subset X$ to $X_i$.
  \end{itemize}
\end{defn}

The concept of metric approximation is needed to define
$L^p$-topology in Section \ref{ssec:Lp}.

\begin{defn}[Compatibility]
  For an asymptotic relation and a metric approximation $\{f_i\}$
  between $\{X_i\}$ and $X$, we consider the following
  \emph{compatibility condition}, (C), between them.
  \begin{itemize}
  \item[(C)] $f_i(x) \to x$ in $\mathcal{X}$ for any $x \in \bigcup_i
    \Dom(f_i)$.
  \end{itemize}
  If the compatibility condition holds, we say that the asymptotic
  relation and the metric approximation $\{f_i\}$ are \emph{compatible
    to each other}.
\end{defn}

\begin{lem} \label{lem:asymprel}
  \begin{enumerate}
  \item For a given (linear) asymptotic relation between $\{X_i\}$ and
    $X$, there exists a (linear) metric approximation $\{f_i\}$
    compatible with it such that $\Dom(f_i) = X$ for all $i$ in the
    nonlinear case and that $\Dom(f_i)$ consists of finite linear
    combinations of a given complete basis on $X$ in the linear case.
  \item For a given (linear) metric approximation $\{f_i\}$ for
    $\{X_i\}$ and $X$, there exists a unique (linear) asymptotic
    relation between $\{X_i\}$ and $X$ compatible with $\{f_i\}$
  \end{enumerate}
\end{lem}

\begin{proof}
  (1): Let an asymptotic relation between $\{X_i\}$ and $X$ be given.
  By (A2), for each $x \in X$ we can choose a net $x_i \in X_i$ such
  that $x_i \to x$ with respect to the asymptotic relation and set
  $f_i(x) := x_i$.  This defines a map $f_i : \Dom(f_i) = X \to X_i$.
  (A3) implies (B2).
  
  Assume that the asymptotic relation is linear and let $\mathcal{B}$
  be a complete linear basis on $X$.  We define a map $f_i$ on
  $\mathcal{B}$ in the same way as above.  Then it extends to a linear
  map from the set of finite linear combinations of $\mathcal{B}$, say
  $\Dom(f_i)$, to $X_i$.  $\Dom(f_i)$ is a dense linear subspace of
  $X$.
  
  (2): Let $\{f_i\}$ be a metric approximation for $\{X_i\}$ and $X$.
  We define a convergence $X_i \ni x_i \to x \in X$ by the following
  condition: There exists a sequence $\tilde x_j \in \bigcup_i
  \Dom(f_i)$, $j=1,2,\dots$, converging to $x$ in $X$ such that
  \[
  \quad \lim_j \limsup_i d_{X_i}(x_i,f_i(\tilde x_j)) = 0
  \]
  (compare \eqref{eq:Lpdef}).  This convergence together with (A1)
  induces a unique topology on $\mathcal{X}$.  We shall show that this
  topology is an asymptotic relation.  In fact, (B1) implies (A2).
  Let us verify (A3).  Assume $X_i \ni x_i \to x \in X$ and $X_i \ni
  y_i \to y \in X$.  Then there exist $\tilde x_j, \tilde y_j \in
  \bigcup_i \Dom(f_i)$ such that
  \begin{align*}
    &\tilde x_j \to x, \qquad
    \lim_j \limsup_i d_{X_i}(x_i,f_i(\tilde x_j)) = 0, \\
    &\tilde y_j \to y, \qquad
    \lim_j \limsup_i d_{X_i}(y_i,f_i(\tilde y_j)) = 0.
  \end{align*}
  Since $\{f_i\}$ is a metric approximation, we have
  \[
  \lim_j \lim_i d_{X_i}(f_i(\tilde x_j),f_i(\tilde y_j))
  = \lim_j d_X(\tilde x_j,\tilde y_j) = d_X(x,y).
  \]
  Therefore the triangle inequalities show that $d_{X_i}(x_i,y_i) \to
  d_X(x,y)$.  (A4) is obtained from the definition of the convergence
  and a triangle inequality.  Thus the topology on $\mathcal{X}$
  defined here is an asymptotic relation.  It is obvious that the
  compatibility condition is satisfied.

  Supposing that $\{f_i\}$ is linear, it is easy to see that
  the asymptotic relation defined above is linear.
\end{proof}

\begin{defn}[Asymptotic continuity]
  Let $\{X_i\}$ and $X$ have an asymptotic relation.  A compatible
  metric approximation $\{f_i\}$ is said to be \emph{asymptotically
    continuous} if $\Dom(f_i) = X$ for any $i$ and if
  $f_i(x_i) \to x$ in $\mathcal{X}$ holds for any net
  $X \ni x_i \to x \in X$.
\end{defn}

\begin{lem} \label{lem:asympconti}
  Assume that $\{X_i\}$ and $X$ have a (linear) asymptotic relation
  and that $X$ is separable.  Then, there exists an asymptotically
  continuous (linear) metric approximation $\{f_i : \Dom(f_i) = X \to
  X_i\}$ compatible with it such that each $f_i$ is a Borel map.
\end{lem}

\begin{proof}
  Let us first consider the nonlinear case.  Take a dense countable
  subset $\{a_k\}_{k \in \N}$ of $X$ and set $A_n :=
  \{a_1,a_2,\dots,a_n\}$.  For a given $x \in X$ we choose $a_k$ with
  smallest $k$ among the points in $A_n$ nearest to $x$; then we set
  $\pi_n(x) := a_k$.  This defines a Borel map $\pi_n : X \to A_n$.
  By Lemma \ref{lem:asymprel}, there exists a metric approximation
  $\{g_i : X \to X_i\}$ compatible with the given asymptotic relation
  between $\{X_i\}$ and $X$.  Let
  \[
  \epsilon_{n,i} := \sup
  \{\; | d_{X_j}(g_j(a),g_j(a')) - d_X(a,a') |\ \mid
  \ j\ge i,\  a,a' \in A_n
  \;\}.
  \]
  Then, for each $n \in \N$ we have $\lim_i \epsilon_{n,i} = 0$.
  Hence, there exists a net $n(i) \to \infty$ such that $\lim_i
  \epsilon_{n(i),i} = 0$.  Define $f_i := g_i \circ \pi_{n(i)}$.
  Since the image of each $\pi_{n(i)}$ consists of finitely many
  points, it is a Borel map.  We take any $x,x'\in X$ and fix them.
  It then follows that
  \begin{align*}
      |\;d_{X_i}(f_i(x),f_i(x')) & - d_X(x,x')\;| \\
      &\le |\;d_{X_i}(f_i(x),f_i(x')) -
      d_X(\pi_{n(i)}(x),\pi_{n(i)}(x'))\;| \\
      &\quad + |\;d_X(\pi_{n(i)}(x),\pi_{n(i)}(x')) - d_X(x,x')\;| \\
      &\le \epsilon_{n(i),i} + d_X(\pi_{n(i)}(x),x) + d_X(\pi_{n(i)}(x'),x')\\
      &\to 0,
  \end{align*}
  so that $\{f_i\}$ is a metric approximation.  The rest of the proof
  is to show the asymptotic continuity of $\{f_i\}$.  Take points
  $x_i,x \in X$ such that $x_i \to x$.  Since $\pi_{n(i)}(x_i)$
  and $\pi_{n(i)}(x)$ both tend to $x$, we have
  $d_X(\pi_{n(i)}(x_i),\pi_{n(i)}(x)) \to 0$, which together with
  $\epsilon_{n(i),i} \to 0$ implies $d_X(f_i(x_i),f_i(x)) \to 0$.  By
  (A4), we obtain $f_i(x_i) \to x$ in $\mathcal{X}$.  This completes
  the proof in the nonlinear case.
  
  Assume that the asymptotic relation is linear.  Since $X$ is
  separable, there is a countable complete linear basis $\mathcal{B} =
  \{a_1,a_2,\dots\}$ of $X$ such that $d_X(o,a_k)$ is bounded away
  from zero and from infinity as $k \to \infty$.  We have a unique
  Hilbert metric on $X$ for which $\mathcal{B}$ is a complete
  orthonormal basis.  It follows that the topology of the Hilbert
  metric on $X$ coincides with that of $d_X$.  Let $A_n$ be the linear
  subspace spanned by $a_1,a_2,\dots,a_n$ and $\pi_n : X \to A_n$ the
  orthogonal projection.  Every $\pi_n$ is continuous and satisfies
  that $\pi_m = \pi_m \circ \pi_n$ for $n \le m$.  By Lemma
  \ref{lem:asymprel}, we find a linear metric approximation $\{g_i\}$
  compatible with the asymptotic relation such that $\Dom(g_i)$
  consists of finite linear combinations of $\mathcal{B}$.  Let
  \[
  \epsilon_{n,i} := \sup
  \{\; | d_{X_j}(g_j(a),g_j(a')) - d_X(a,a') |\ \mid
  \ j\ge i,\  a,a' \in A_n \cap B(o,1/n)
  \;\}.
  \]
  where $o \in X$ is the origin.  The rest of the proof is same as
  in the nonlinear case above.
\end{proof}

\begin{lem} \label{lem:asympconti2}
  Assume that $\{X_i\}$ and $X$ have an asymptotic relation and
  $\{f_i\}$ an asymptotically continuous metric approximation
  compatible with it.  Then, for any compact subset $C \subset X$ we
  have $\dis f_i|_C \to 0$.
\end{lem}

\begin{proof}
  If not, there exist a number $\delta > 0$ and nets $x_i,y_i \in C$
  such that
  \begin{equation}
    \label{eq:asympcontiC}
    |d_{X_i}(f_i(x_i),f_i(y_i)) - d_X(x_i,y_i)| \ge \delta.
  \end{equation}
  Since $C$ is compact, by replacing the nets by subnets, we may assume
  that $x_i \to x$, $y_i \to y$ for some points $x,y \in C$.
  Then we have $d_X(x_i,y_i) \to d_X(x,y)$ and, by the
  asymptotic continuity of $\{f_i\}$,
  \[
  d_{X_i}(f_i(x_i),f_i(x)) \to 0, \qquad d_{X_i}(f_i(y_i),f_i(y)) \to
  0.
  \]
  Moreover, $d_{X_i}(f_i(x),f_i(y)) \to
  d_X(x,y)$.  Thus, by using triangle inequalities, the left-hand
  side of \eqref{eq:asympcontiC} tends to zero.  This is a
  contradiction.
\end{proof}

We shall briefly mention a connection between asymptotic relation and
Gromov-Hausdorff convergence.

\begin{defn}[Asymptotic compactness]
  Assume that $\{X_i\}$ and $X$ have an asymptotic relation.  We say
  that a net $x_i \in X_i$ is \emph{bounded} if $d_{X_i}(x_i,o_i)$ is
  bounded for some convergent net $o_i \in X_i$.  The asymptotic
  relation is said to be \emph{asymptotically compact} if any bounded
  net $x_i \in X_i$ has a convergent subnet in $\mathcal{X}$ with
  respect to the asymptotic relation.
\end{defn}

\begin{prop} \label{prop:proper}
  If there exists an asymptotically compact asymptotic relation
  between $\{X_i\}$ and $X$, then $X$ is proper.
\end{prop}

\begin{prop} \label{prop:GHasympcpt}
  Let $X_i$ and $X$ be proper metric spaces, let $o_i \in X_i$ be a
  net, and let $o \in X$ be a point.  Then, the pointed space
  $(X_i,o_i)$ converges to $(X,o)$ in the Gromov-Hausdorff topology
  iff there exists an asymptotically compact asymptotic relation
  between $\{X_i\}$ and $X$ for which $o_i$ converges to $o$.
\end{prop}

The proofs of Propositions \ref{prop:proper} and \ref{prop:GHasympcpt}
are easy.  Also, they both follow from Theorem
\ref{thm:cpt} blow by setting $E_i := 0$ and $E := 0$.

\subsection{Asymptotic relation between $L^p$-spaces} \label{ssec:Lp}

Let $M$ be a measure space (as we said in Section \ref{sec:prelim}, it
is a locally compact Polish space with a full supported Radon
measure), $Y$ a metric space, and $p \ge 1$ a real number.  Given two
measurable maps $u,v : M \to Y$, we define the \emph{$L^p$-distance
  $d_{L^p}(u,v)$} between them by
\[
d_{L^p}(u,v) := \left(\int_M d_Y(u(x),v(x))^p \; dx\right)^{\frac1p}
\le +\infty,
\]
where $\int_M dx$ means the integrating over $M$ by the measure on
$M$.  For a measurable map $\xi : M \to Y$ we define
\[
L^p_\xi(M,Y) := \{\; u : M \to Y
\mid \text{measurable map with $d_{L^p}(u,\xi) < +\infty$} \;\}.
\]
We identify two maps in $L^p_\xi(M,Y)$ if they are equal a.e. on $M$,
so that $L^p_\xi(M,Y)$ becomes a metric space with metric $d_{L^p}$.
If $Y$ is complete (resp.~separable), then $L^p_\xi(M,Y)$ is also
complete (resp.~separable).

We first introduce a natural asymptotic relation between $L^p$-spaces
in the case where the target space is a Banach space.  Let
$\mathbb{B}$ be a Banach space over $\R$ with norm $\|\cdot\|$ and
origin $o$.  Note that $L^p_o(M,\mathbb{B})$ is a Banach space with
respect to $d_{L^p}$.  The \emph{support} `$\supp u$' of a measurable
map $u : M \to Y$ is defined to be the subset of $M$ satisfying the
condition that $x \in M \setminus \supp u$ iff there exists an open
neighborhood $U$ of $x$ such that $u = o$ a.e.~on $U$.  Denote by
$C_o(M,\mathbb{B})$ the set of continuous maps $u : M \to \mathbb{B}$
with compact support $\supp u$.  We have the following lemma in a
standard way and the proof is omitted.

\begin{lem} \label{lem:dense}
  $C_o(M,\mathbb{B})$ is a dense linear subspace of $L^p_o(M,\mathbb{B})$.
\end{lem}

Let $\{\varphi_i : M_i \supset \Dom(\varphi_i) \to M\}$ be a measure
approximation for measure spaces $M_i$ and $M$.
For $u \in C_o(M,\mathbb{B})$, we define
\begin{equation}
  \label{eq:Phii}
  \Phi_iu(x) := 
  \begin{cases}
    u \circ \varphi_i(x) &\text{for $x \in \Dom(\varphi_i)$,}\\
    o &\text{for $x \in M_i \setminus \Dom(\varphi_i)$}.
  \end{cases}
\end{equation}
For each $i$, $\Phi_i u : M_i \to \mathbb{B}$ is a measurable map and
$\Phi_i$ is a linear map from $C_o(M,\mathbb{B})$.  Define
$\Dom(\Phi_i)$ to be the set of $u \in C_o(M,\mathbb{B})$ such that
$d_{L^p}(\Phi_ju,o) < +\infty$ for any $j \ge i$.  It is easy to prove
that $\Dom(\Phi_i)$ is a linear subspace of $C_o(M,\mathbb{B})$.

\begin{prop} \label{prop:metapprox}
  We have
  \begin{enumerate}
  \item $\bigcup_i \Dom(\Phi_i) = C_o(M,\mathbb{B})$;
  \item $\{\Phi_i : L^p_o(M,\mathbb{B}) \supset \Dom(\Phi_i) \to
    L^p_o(M_i,\mathbb{B})\}$ is a linear metric approximation.
  \end{enumerate}
\end{prop}

\begin{proof}
  We take any two maps $u,v \in C_o(M,\mathbb{B})$ and fix them.
  Since $M \ni x \mapsto \|u(x) - v(x)\|^p$ is a continuous function
  with compact support, it follows from Definition
  \ref{defn:mapprox}(M2) that
  \begin{align*}
    d_{L^p}(\Phi_iu,\Phi_iv)^p &= \int_{\Dom(\varphi_i)} \|u \circ
    \varphi_i(x) - v \circ \varphi_i(x)\|^p \; dx \\
    &\to \int_M \|u(x) - v(x)\|^p\;dx = d_{L^p}(u,v)^p.
  \end{align*}
  In particular,
  \[
  \lim_i d_{L^p}(\Phi_iu,o)
  = \lim_i d_{L^p}(\Phi_iu,\Phi_io)
  =  d_{L^p}(u,o) < +\infty,
  \]
  which implies (1).  By recalling that $C_o(M,\mathbb{B})$ is dense
  in $L^p_o(M,\mathbb{B})$, we obtain (2).  This completes the proof.
\end{proof}

\begin{defn}[$L^p$-topology] \label{defn:LptopB}
  As is seen in Lemma \ref{lem:asymprel}, the metric approximation
  $\{\Phi_i\}$ induces a unique linear asymptotic relation
  between $X_i := L^p_o(M_i,\mathbb{B})$ and $X :=
  L^p_o(M,\mathbb{B})$.  We call the topology on $\mathcal{X} =
  (\bigsqcup_i X_i) \sqcup X$ the \emph{$L^p$-topology} and a
  convergence for it an \emph{$L^p$-convergence}.
\end{defn}

\begin{rem} \label{rem:Lptop}
  In the proof of Proposition \ref{prop:metapprox}, the continuity of
  the maps $u$ and $v$ is necessary, so that the dense property of
  $C_o(M,\mathbb{B})$ in $L^p_o(M,\mathbb{B})$ is important to define
  the $L^p$-topology.  If the target $Y$ is a general metric space,
  then $C_o(M,Y)$ is not necessarily dense in $L^p_o(M,Y)$ (e.g. in
  the case where $Y$ is disconnected).
\end{rem}

We next consider the case where the targets are spaces of a
Gromov-Hausdorff convergent net.  Let $(Y_i,y_i)$ and $(Y,y_0)$ be
pointed proper metric spaces such that $(Y_i,y_i)$ converges to
$(Y,y_0)$ in the pointed Gromov-Hausdorff topology.

\begin{lem}
  There exists a separable real Banach space $(\mathbb{B},\|\cdot\|)$
  into which all $Y_i$ and $Y$ are embedded isometrically.
\end{lem}

\begin{proof}
  By Proposition \ref{prop:GHunion}, all $Y_i$ and $Y$ can
  isometrically be embedded into the metric space $\mathcal{Y} :=
  (\bigsqcup_i Y_i) \sqcup Y$ with metric $d_{\mathcal{Y}}$.  As is
  well-known, we can isometrically embed $\mathcal{Y}$ into the real
  Banach space, say $(\mathbb{B},\|\cdot\|)$, consisting of continuous
  bounded functions on $\mathcal{Y}$ with uniform norm.  The embedding
  map is
  \[
  \mathcal{Y} \ni y \mapsto
  d_{\mathcal{Y}}(y,\cdot) - d_{\mathcal{Y}}(y_0,\cdot) \in \mathbb{B}.
  \]
  Since $\mathcal{Y}$ is separable, so is $\mathbb{B}$.    
\end{proof}

We take $(\mathbb{B},\|\cdot\|)$ as in the lemma.
Let $\xi_i : M_i \to Y_i \subset \mathbb{B}$ be measurable maps and
$\xi : M \to Y \subset \mathbb{B}$ a continuous map such that
\begin{equation}
  \label{eq:xi}
  \lim_i \int_{\Dom(\varphi_i)} \|\xi_i(x) - \xi\circ\varphi_i(x)\|^p
  \; dx = 0,
\end{equation}
and consider $L^p_{\xi_i}(M_i,Y_i)$ and
$L^p_\xi(M,Y)$.  We embed them to $L^p_o(M_i,\mathbb{B})$ and
$L^p_o(M,\mathbb{B})$ respectively by
\begin{equation}
  \label{eq:embLp}
  \begin{split}
    \iota_i : L^p_{\xi_i}(M_i,Y_i) \ni u
    &\mapsto u-\xi_i \in L^p_o(M_i,\mathbb{B}),\\
    \iota : L^p_\xi(M,Y) \ni u &\mapsto u-\xi \in L^p_o(M,\mathbb{B}),
  \end{split}
\end{equation}
both which are isometric.

\begin{rem} \label{rem:linear}
  If $Y_i$ and $Y$ are real Banach spaces (not necessarily isometric
  to each other) and if $\xi_i$ and $\xi$ are their origins, then the
  embeddings \eqref{eq:embLp} are linear maps.
\end{rem}

\begin{defn}[$L^p$-topology] \label{defn:Lptop}
  We define the \emph{$L^p$-topology on}
  \[
  \left( \bigsqcup_i L^p_{\xi_i}(M_i,Y_i) \right) \sqcup L^p_\xi(M,Y)
  \]
  as the restriction of that of $(\bigsqcup_i L^p_o(M_i,\mathbb{B}))
  \sqcup L^p_o(M,\mathbb{B})$.
\end{defn}

In the following, let us see that the $L^p$-topology defined here has
some natural properties.

We need a lemma.  $|A|$ denotes the measure of a measurable set $A$.

\begin{lem} \label{lem:OF}
  \begin{enumerate}
  \item For any open subset $O \subset M$ we have
    \[
    \liminf_i |\varphi_i^{-1}(O)| \ge |O|.
    \]
  \item For any closed subset $F \subset M$ we have
    \[
    \limsup_i |\varphi_i^{-1}(F)| \le |F|.
    \]
  \end{enumerate}
\end{lem}

We omit the proof of Lemma \ref{lem:OF}.


We give two measure approximations $\{\varphi_i : M_i \supset
\Dom(\varphi_i) \to M\}$ and $\{\psi_i : M_i \supset \Dom(\psi_i) \to
M\}$.

\begin{defn}[Equivalence relation between measure approximations]
  We say that $\{\varphi_i\}$ and $\{\psi_i\}$ are \emph{equivalent}
  if for any compact subset $C \subset M$ and for any $\epsilon > 0$
  there exists $i(C,\epsilon)$ such that for any $i \ge i(C,\epsilon)$
  and any $x \in \Dom(\varphi_i) \cap \Dom(\psi_i) \cap (\varphi_i^{-1}(C)
  \cup \psi_i^{-1}(C))$ we have $d_M(\varphi_i(x),\psi_i(x)) <
  \epsilon$, where $d_M$ is a distance function on $M$ compatible with
  the topology of $M$.
\end{defn}

Note that the equivalence relation defined here is independent of
the distance function $d_M$, because it can be described only by
the uniform structure on $M$.

\begin{lem} \label{lem:measapproxequiv}
  If $\{\varphi_i\}$ and $\{\psi_i\}$ are equivalent, then the two
  $L^p$-topologies induced from $\{\varphi_i\}$ and $\{\psi_i\}$
  coincide.
\end{lem}

\begin{proof}
  Assume that $\{\varphi_i\}$ and $\{\psi_i\}$ are equivalent.  Take
  any $u \in C_o(M,\mathbb{B})$ and fix it.  We set $C := \supp u$.
  By the uniform continuity of $u$, for any $\epsilon > 0$ there is
  $i(C,\epsilon)$ such that for any $x \in \Dom(\varphi_i) \cap
  \Dom(\psi_i) \cap (\varphi_i^{-1}(C) \cup \psi_i^{-1}(C))$, $i \ge
  i(C,\epsilon)$, we have $\|u\circ\varphi_i(x)-u\circ\psi_i(x)\| <
  \epsilon$ and therefore,
  \[
  d_{L^p}(\Phi_i(u),\Psi_i(u))^p
  \le \epsilon^p \, |\varphi_i^{-1}(C) \cup \psi_i^{-1}(C)|,
  \]
  where $\Psi_i$ is defined for $\psi_i$ in the same manner as
  \eqref{eq:Phii}.  By using Lemma \ref{lem:OF}(2), this implies that
  $\limsup_i d_{L^p}(\Phi_i(u),\Psi_i(u)) \le 2\epsilon^p |C|$.  By
  the arbitrariness of $\epsilon$ we have $\lim_i
  d_{L^p}(\Phi_i(u),\Psi_i(u)) = 0$.  Since this holds for any $u \in
  C_o(M,\mathbb{B})$, the two measure approximations $\{\Phi_i\}$ and
  $\{\Psi_i\}$ induce the same asymptotic relation between
  $L^p_o(M_i,\mathbb{B})$ and $L^p_o(M,\mathbb{B})$.  The restriction
  in Definition \ref{defn:Lptop} has the same topology.  This
  completes the proof.
\end{proof}

\begin{lem} \label{lem:indepB}
  The $L^p$-topology of Definition \ref{defn:Lptop} is independent of
  the Banach space $\mathbb{B}$ and the embedding $Y_i, Y
  \hookrightarrow \mathbb{B}$.
\end{lem}

\begin{proof}
  Assume we have two embeddings $\alpha : Y_i, Y \hookrightarrow
  \mathbb{B}_\alpha$ and $\beta : Y_i, Y \hookrightarrow
  \mathbb{B}_\beta$ into two real Banach spaces.  Let $Z$ be the gluing of
  $\mathbb{B}_\alpha$ and $\mathbb{B}_\beta$ along the isometric
  images $\alpha(Y)$ and
  $\beta(Y)$, and $d_Z$ the metric on $Z$ defined by
  \[
  d_Z(x,y) :=
  \begin{cases}
    \|x-y\| &\text{if $x,y \in \mathbb{B}_\alpha$
      or if $x,y \in \mathbb{B}_\beta$,}\\
    \inf_{z \in \mathcal{Y}}
    (\|x-z\| + \|z-y\|) &\text{otherwise,}
  \end{cases}
  \]
  for $x,y \in Z$.  Then we embed $(Z,d_Z)$ into a real Banach space
  $\mathbb{B}$ isometrically, so that we have the embeddings $Y_i,Y
  \hookrightarrow \mathbb{B}_\alpha, \mathbb{B}_\beta \hookrightarrow
  \mathbb{B}$.  Identifying $\mathbb{B}_\alpha$ and $\mathbb{B}_\beta$
  with their isometric images in $\mathbb{B}$, we assume that $\alpha$
  and $\beta$ each take their values in $\mathbb{B}$.  Let $u_i \in
  L^p_{\xi_i}(M_i,Y_i)$ and set
  \begin{gather*}
    \hat u_i^\alpha := \alpha \circ u_i,
    \ \hat u_i^\beta := \beta \circ u_i
    \ \in L^p_{\xi_i}(M_i,\mathbb{B}),\\
    u_i^\alpha := \alpha \circ u_i - \alpha \circ \xi_i, 
    \ u_i^\beta := \beta \circ u_i - \beta \circ \xi_i
    \ \in L^p_o(M_i,\mathbb{B}).
  \end{gather*}
  Take a map $u \in C_o(M,\mathbb{B})$ with $|\supp u| > 0$, and set
  $\hat u := u + \xi$.  Assume that
  \[
  \limsup_i d_{L^p}(u_i^\alpha,\Phi_i(u)) < \epsilon
  \]
  for a number $\epsilon > 0$.
  Let us estimate $d_{L^p}(u_i^\beta,\Phi_i(u))$ by $\epsilon$.
  We have
  \begin{align*}
      d_{L^p}(u_i^\beta,\Phi_i(u))^p
      &= \int_{\varphi_i^{-1}(\supp u)} \|u_i^\beta(x) -
      u\circ\varphi_i(x)\|^p\,dx\\
      &\quad + \int_{M_i \setminus \varphi_i^{-1}(\supp u)}
      \|u_i^\beta(x)\|^p \,dx.
  \end{align*}
  Since $\|u_i^\beta(x)\| = d_{Y_i}(u_i(x),\xi_i(x)) =
  \|u_i^\alpha(x)\|$, the $\beta$ of the second term of the right
  hand-side of the above can be replaced by $\alpha$ and is less than
  $\epsilon^p$ for $i$ large enough.  We estimate the first term.
  Note that, by \eqref{eq:xi}, the limsup of the first term is less
  than or equal to
  \begin{equation}
    \label{eq:indepB0}
    \limsup_i \int_{\varphi_i^{-1}(\supp u)} \|\hat u_i^\beta(x) -
    \hat u\circ\varphi_i(x)\|^p\,dx.
  \end{equation}
  Let $R$ be a number with $R > \sup_{z \in \hat u(\supp u)}
  \|\alpha(y_0)-z\|$.  Note that $\alpha(y_0) = \beta(y_0)$.  There is
  $i(\epsilon)$ such that $\|\alpha(y)-\beta(y)\| < \epsilon\,|\supp
  u|^{-1/p}$ for any $i \ge i(\epsilon)$ and any $y \in B(y_i,2R)
  (\subset Y_i)$.  In particular, if a point $x \in M_i$ for $i \ge
  i(\epsilon)$ satisfies
  $u_i(x) \in B(y_i,2R)$, then
  \begin{equation}
    \label{eq:indepB1}
    \|\hat u_i^\alpha(x) - \hat u_i^\beta(x)\| <
    \epsilon\,|\supp u|^{-1/p}.
  \end{equation}
  Denote by $S_{\epsilon,i}$ the set of all $x \in
  \varphi_i^{-1}(\supp u)$ such that $u_i(x) \not\in B(y_i,2R)$.
  Then, \eqref{eq:indepB1} holds for any $x \in M_i \setminus
  S_{\epsilon,i}$, $i \ge i(\epsilon)$.  Since $\liminf_i \inf_{x \in
    S_{\epsilon,i}} \|\hat u_i^\alpha(x)-\alpha(y_0)\| \ge 2R$ and
  $\hat u\circ\varphi_i(x) \in \hat u(\supp u)$, the triangle
  inequality implies that $\liminf_i \inf_{x \in S_{\epsilon,i}}
  \|\hat u_i^\alpha(x)-\hat u\circ\varphi_i(x)\| \ge R$ and so
  \begin{equation}
    \label{eq:indepB2}
    R^p \limsup_i |S_{\epsilon,i}| \le 
    \limsup_i \int_{S_{\epsilon,i}}
    \|\hat u_i^\alpha(x)-\hat u\circ\varphi_i(x)\|^p \, dx
    \le \limsup_i d_{L^p}(u_i^\alpha,\Phi_i(u))^p < \epsilon^p.
  \end{equation}
  Since $\|\hat u_i^\alpha(x) - \alpha(y_i)\| = \|\hat u_i^\beta(x) -
  \beta(y_i)\|$,
  we have
  \begin{equation}
    \label{eq:indepB3}
    \limsup_i \bigl|\;\|\hat u_i^\alpha(x)- \hat u\circ\varphi_i(x)\|
    -\|\hat u_i^\beta(x)- \hat u\circ\varphi_i(x)\|\;\bigr| \le 2R.
  \end{equation}
  By \eqref{eq:indepB1}, \eqref{eq:indepB2}, and \eqref{eq:indepB3},
  we see that \eqref{eq:indepB0} is
  \begin{align*}
    &\le 
    \limsup_i \Biggl\{\int_{\varphi_i^{-1}(\supp u) \setminus S_{\epsilon,i}}
    (\|\hat u_i^\alpha(x) -
    \hat u\circ\varphi_i(x)\| + \epsilon\,|\supp u|^{-1/p})^p\,dx\\
    &\quad + \int_{\varphi_i^{-1}(\supp u) \cap S_{\epsilon,i}}
    (\|\hat u_i^\alpha(x) - \hat u\circ\varphi_i(x)\| + 2R)^p\,dx
    \Biggr\}\\
    &\le \limsup_i \int_{\varphi_i^{-1}(\supp u)}
    2^{p-1}\|\hat u_i^\alpha(x) - \hat u\circ\varphi_i(x)\|^p \, dx
    + 2^{p-1}(\epsilon^p + 2^p\epsilon^p)\\
    &\le (2^p+2^{2p-1})\epsilon^p.
  \end{align*}
  Since $\epsilon$ is arbitrary, we obtain
  \[
  \limsup_i d_{L^p}(u_i^\beta,\Phi_i(u))
  \le C_p \limsup_i d_{L^p}(u_i^\alpha,\Phi_i(u)),
  \]
  where $C_p$ is a constant depending only on $p$.
  
  Now, assume that $u_i^\alpha \in L^p_o(M_i,\mathbb{B})$
  $L^p$-converges to a map $u \in L^p_o(M,\mathbb{B})$.  There is a
  continuous approximation $\tilde u_j \in C_o(M,\mathbb{B})$ of $u$,
  $j = 1,2,\dots$, such that $|\supp\tilde u_j| > 0$ and
  $d_{L^p}(\tilde u_j,u) \to 0$ as $j \to \infty$.  Since
  $u_i^\alpha \to u$ in $L^p$, we have
  \[
  \lim_j \limsup_i d_{L^p}(u_i^\beta,\Phi_i(\tilde u_j)) \le
  C_p \lim_j \limsup_i d_{L^p}(u_i^\alpha,\Phi_i(\tilde u_j)) = 0,
  \]
  which implies that $u_i^\beta$ $L^p$-converges to $u$.
  This completes the proof.
\end{proof}

\begin{rem} \label{rem:xi}
  The condition \eqref{eq:xi} depends on the embedding $Y_i, Y
  \hookrightarrow \mathbb{B}$.  Assume \eqref{eq:xi} holds for some
  embedding $Y_i, Y \hookrightarrow \mathbb{B}$, and take another
  embedding $Y_i, Y \hookrightarrow \mathbb{B}'$.  Then, replacing
  each $\Dom(\varphi_i)$ with a suitable subset, we can keep
  \eqref{eq:xi} for $Y_i, Y \hookrightarrow \mathbb{B}'$.  Since the
  restriction of $\varphi_i$ is equivalent to the original one, and by
  Lemma \ref{lem:measapproxequiv}, these induce the same
  $L^p$-topology.
\end{rem}

Fix an embedding $Y_i,Y \hookrightarrow \mathbb{B}$.
Let us prove that $(\bigsqcup_i L^p_{\xi_i}(M_i,Y_i)) \sqcup
L^p_\xi(M,Y)$ is closed in $(\bigsqcup_i L^p_o(M_i,\mathbb{B})) \sqcup
L^p_o(M,\mathbb{B})$.

\begin{lem} \label{lem:closed}
  Let $u_i \in L_{\xi_i}^p(M_i,Y_i)$ be a net such that $v_i :=
  u_i - \xi_i$ converges to a map $v \in
  L_o^p(M,\mathbb{B})$ in the $L^p$-topology.  Then, $u :=
  v + \xi$ belongs to $L_\xi^p(M,Y)$.
\end{lem}

\begin{proof}
  Suppose that $u = v + \xi$ does not belongs to $L_\xi^p(M,Y)$.
  Then, $|\{ x \in M \mid u(x) \not\in Y\}| > 0$.  There is a small
  number $\delta > 0$ such that $U := \{ x \in M \mid
  d_{\mathbb{B}}(u(x),Y) > \delta \}$ has positive measure, where
  $d_{\mathbb{B}}$ is the distance function induced from the Banach
  norm $\|\cdot\|$ on $\mathbb{B}$.  For the $v$ and any $\epsilon >
  0$, there exists a continuous map $\tilde v_\epsilon \in
  C_o(M,\mathbb{B})$ such that $d_{L^p}(v,\tilde v_\epsilon) <
  \epsilon$.  Setting $\tilde u_\epsilon := \tilde v_\epsilon + \xi$
  we have
  \begin{align*}
    \epsilon^p &> \int_{d_{\mathbb{B}}(u(x),Y) > \delta,\ d_{\mathbb{B}}(\tilde u_\epsilon(x),Y)
      \le \delta/2}
    \| u(x) - \tilde u_\epsilon(x)\|^p \; dx\\
    &\ge (\delta/2)^p |\{x \in U \mid d_{\mathbb{B}}(\tilde u_\epsilon(x),Y)
    \le \delta/2 \}|
  \end{align*}
  and therefore
  \[
  \liminf_{\epsilon \to 0+} |\{ x \in M \mid d_{\mathbb{B}}(\tilde u_\epsilon(x),Y)
  > \delta/2 \}| > 0.
  \]
  Since $\tilde u_\epsilon$ is continuous, applying Lemma \ref{lem:OF}(1)
  yields
  \[
  \liminf_{\epsilon \to 0+} \liminf_i
  |\{ x \in \Dom(\varphi_i) \mid d_{\mathbb{B}}(\tilde u_\epsilon \circ \varphi_i(x),Y)
  > \delta/2\}| > 0,
  \]
  which implies
  \[
  \liminf_{\epsilon\to 0+} \liminf_i \int_{\Dom(\varphi_i)} \|u_i(x) -
  \tilde u_\epsilon \circ \varphi_i(x)\|^p \; dx > 0.
  \]
  This is a contradiction to $v_i \to v$ in the
  $L^p$-topology and \eqref{eq:xi}.
\end{proof}

\begin{lem} \label{lem:psi}
  Let $(Z,d_Z)$ be a separable metric space.  For any closed subset $F
  \subset Z$ there exists a Borel map $\psi : Z \to F$ such that
  $d_Z(\psi(z),z) \le 2 d_Z(z,F)$ for any $z \in Z$.
\end{lem}

\begin{proof}
  It is necessary that $\psi(z) = z$ for $z \in F$.  We shall define
  $\psi$ on $Z \setminus F$.  Take a dense countable subset
  $\{z_n\}_{n \in \N} \subset Z \setminus F$ and set $r_n :=
  d_Z(z_n,F)$ and $B_n := B(z_n,r_n/2)$.  $\{B_n\}_{n \in \N}$ is an
  open covering of $Z \setminus F$.  Let $B_n' := B_n \setminus
  \bigcup_{k=1}^{n-1} B_k$.  Then, $\{B_n'\}_{n\in\N}$ is a covering
  of $Z \setminus F$ consisting of disjoint Borel sets.  For each $n
  \in \N$ there is a point $a_n \in F$ such that $d_Z(a_n,z_n) <
  3r_n/2$.  We set $\psi(z) := a_n$ for $z \in B_n'$.  This defines
  a Borel map $\psi : Z \to F$.  The triangle inequality shows that
  $d_Z(\psi(z),z) \le 2 d_Z(z,F)$ for any $z \in Z$.
\end{proof}

The following lemma verifies (A2) of Definition \ref{defn:asymprel}.

\begin{lem} \label{lem:a2}
  For any $u \in L^p_\xi(M,Y)$ there exists a net $u_i \in
  L^p_{\xi_i}(M_i,Y_i)$ converging to $u$ in the $L^p$-topology.
\end{lem}

\begin{proof}
  Let $u \in L^p_\xi(M,Y)$.  By Lemma \ref{lem:dense}, for any
  $\epsilon > 0$ there exists a map $\tilde u_\epsilon \in
  C(M,\mathbb{B})$ such that $\tilde u_\epsilon - \xi \in
  C_o(M,\mathbb{B})$ and $d_{L^p}(\tilde u_\epsilon,u)^p < \epsilon$.
  Set $C_\epsilon := \supp(\tilde u_\epsilon - \xi)$.
  Note that $C_\epsilon$ and $\tilde u_\epsilon(C_\epsilon)$ are compact.
  Applying Lemma \ref{lem:psi} we have a Borel map
  $\psi_i : \mathbb{B} \to Y_i$ such that
  $\|\psi(y)-y\| \le 2 d_{\mathbb{B}}(y,Y_i)$ for any $y \in \mathbb{B}$.
  We define, for $x \in M_i$,
  \[
  u_{\epsilon,i}(x) :=
  \begin{cases}
    \psi_i \circ \tilde u_\epsilon \circ \varphi_i(x)
    &\text{if $x \in \varphi_i^{-1}(C_\epsilon)$},\\
    \xi_i(x) &\text{if $x \in M_i \setminus \varphi_i^{-1}(C_\epsilon)$}.
  \end{cases}
  \]
  It then follows that $u_{\epsilon,i}(M_i) \subset Y_i$ and
  \begin{align*}
    &\int_{\Dom(\varphi_i)} \|u_{\epsilon,i}(x)
    - \tilde u_\epsilon \circ \varphi_i(x)\|^p\;dx\\
    &= \int_{\varphi_i^{-1}(C_\epsilon)}
    \|\psi_i \circ \tilde u_\epsilon \circ \varphi_i(x)
    - \tilde u_\epsilon \circ \varphi_i(x)\|^p\;dx\\
    &\ \ + \int_{\Dom(\varphi_i) \setminus \varphi_i^{-1}(C_\epsilon)}
    \|\xi_i(x) - \xi\circ\varphi_i(x)\|^p\;dx.
  \end{align*}
  By \eqref{eq:xi}, the second term of the right-hand side tends to
  zero as $i \to \infty$.
  The first term is
  \begin{equation}
    \label{eq:2p}
    \le 2^p\int_{\varphi_i^{-1}(C_\epsilon)}
    d_{\mathbb{B}}(\tilde u_\epsilon(x) \circ \varphi_i(x),Y_i)^p\;dx.
  \end{equation}
  Since $\tilde u_\epsilon(C_\epsilon)$ is compact and $(Y_i,y_i) \to (Y,y_0)$,
  we have
  \[
  \lim_i \sup_{y \in \tilde u_\epsilon(C_\epsilon)}
  |d_{\mathbb{B}}(y,Y_i)-d_{\mathbb{B}}(y,Y)| = 0.
  \]
  Therefore, the limit of \eqref{eq:2p} is
  \begin{align*}
    &= \lim_i 2^p\int_{\varphi_i^{-1}(C_\epsilon)}
    d_{\mathbb{B}}(\tilde u_\epsilon(x) \circ \varphi_i(x),Y)^p\;dx\\
    &= 2^p\int_{C_\epsilon}
    d_{\mathbb{B}}(\tilde u_\epsilon(x),Y)^p\;dx
    \le 2^p \,d_{L^p}(\tilde u_\epsilon,u)^p < 2^p \epsilon.
  \end{align*}
  Thus, setting $\hat\Phi_iu := \Phi_i(u-\xi)+\xi_i$, we obtain
  \[
  \lim_i d_{L^p}(u_{\epsilon,i},\hat\Phi_i\tilde u_\epsilon)^p
  =\lim_i \int_{\Dom(\varphi_i)} \|u_{\epsilon,i}(x)
  - \tilde u_\epsilon \circ \varphi_i(x)\|^p\;dx < 2^p\,\epsilon,
  \]
  so that there exists a net of positive numbers $\epsilon_i \to 0$
  such that 
  \[
  d_{L^p}(u_{\epsilon_i,i},\hat\Phi_i\tilde u_{\epsilon_i})^p
  < 2^p\,\epsilon_i
  \]
  for any $i$.  Since $\hat\Phi_i\tilde u_{\epsilon_i}$
  converges to $u$ in the $L^p$-topology, so does $u_i :=
  u_{\epsilon_i,i}$.  This completes the proof.
\end{proof}

We summarize Lemmas \ref{lem:indepB}, \ref{lem:closed}, \ref{lem:a2},
and Remarks \ref{rem:xi}, \ref{rem:linear} into the following:

\begin{thm}
  The $L^p$-topology on $(\bigsqcup_i L^p_{\xi_i}(M_i,Y_i)) \sqcup
  L^p_\xi(M,Y)$ defined in Definition \ref{defn:Lptop} is
  \begin{enumerate}
  \item independent of a representative $\{\varphi_i\}$ of the
    equivalence class of the measure approximations, {\rm(}but may
    depend on an equivalence class{\rm)},
  \item independent of the Banach space $\mathbb{B}$ and the
    embedding $Y_i, Y \hookrightarrow \mathbb{B}$,
  \item an asymptotic relation.
  \end{enumerate}
  Moreover, $(\bigsqcup_i L^p_{\xi_i}(M_i,Y_i)) \sqcup L^p_\xi(M,Y)$
  is closed in $(\bigsqcup_i L^p_o(M_i,\mathbb{B})) \sqcup
  L^p_o(M,\mathbb{B})$.  If $Y_i$ and $Y$ are real Banach spaces and
  if $\xi_i$ and $\xi$ are their origins, then the $L^p$-topology is
  linear as an asymptotic relation.
\end{thm}

For a pointed Gromov-Hausdorff convergent net $M_i \to M$ of proper
pointed measured metric spaces, we have a unique associated
equivalence class (as measure approximations) of measured pointed
Gromov-Hausdorff approximations, so that the $L^p$-topology is
uniquely determined.

We have the $L^p$-topology on the disjoint union
\[
\bigsqcup_{M, (Y,y_0)} L^p_{y_0}(M,Y),
\]
where $M$ runs over all isomorphic classes of
pointed proper measured metric spaces,
and $(Y,y_0)$ runs over all isometric classes of
pointed proper metric spaces.
We can think the union above as the fiber space over
the product space $\{(M,(Y,y_0))\}$ of $\{M\}$ and $\{(Y,y_0)\}$ with
fibers $L^p_{y_0}(M,Y)$.

\section{Variational convergence over metric spaces} \label{sec:varconvmet}

We first give some definitions for variational convergences, e.g.
asymptotic compactness and $\Gamma$-convergence, in Section
\ref{ssec:def}.  In the second section, \ref{ssec:Pcov}, we prove the
asymptotic compactness under a bound of the Poincar\'e constants and a
bound of local covering orders, which is one of the main results of
this paper.  In the third section, \ref{ssec:GHsublevel}, we prove the
equivalence between the compact convergence and the Gromov-Hausdorff
convergence of the energy-sublevel sets.

\subsection{Basics} \label{ssec:def}

Let $X_i$ and $X$ be metric spaces and an asymptotic relation between
them be given.  We take functions $E_i : X_i \to [\,0,+\infty\,]$ and
$E : X \to [\,0,+\infty\,]$.
The following is the Rellich compactness property.

\begin{defn}[Compactness]
  We say that $E$ is \emph{compact} if for any bounded net $x_j \in X$
  with $\limsup E(x_j) < +\infty$ there exists a convergent subnet of
  $\{x_j\}$.
\end{defn}

The next lemma is obvious.

\begin{lem}
  The following (1) and (2) are equivalent.
  \begin{enumerate}
  \item $E$ is lower semi-continuous and compact.
  \item $\{x \in X \mid E(x) \le a\}$ is proper for any $a \ge 0$.
  \end{enumerate}
\end{lem}

The Rellich compactness is generalized into the following.

\begin{defn}[Asymptotic compactness, \cite{Ms:compmedia}]
  The net $\{E_i\}$ of functions is said to be \emph{asymptotically
    compact} if for any bounded net $x_i \in X_i$ with $\limsup
  E_i(x_i) < +\infty$ there exists a convergent subnet of $\{x_i\}$.
\end{defn}

\begin{defn}[$\Gamma$-convergence] \label{defn:Gamma}
  We say that \emph{$E_i$ $\Gamma$-converges to $E$} if the following
  ($\Gamma$1) and ($\Gamma$2) are satisfied:
  \begin{itemize}
  \item[($\Gamma$1)] For any $x \in X$ there exists a net $x_i \in
    X_i$ such that $x_i \to x$ and $E_i(x_i) \to E(x)$.
  \item[($\Gamma$2)] If $X_i \ni x_i \to x \in X$ then $E(x) \le
    \liminf_i E_i(x_i)$.
  \end{itemize}
\end{defn}

\begin{defn}[Compact convergence]
  We say that \emph{$E_i$ compactly converges to $E$} if $E_i$
  $\Gamma$-converges to $E$ and if $\{E_i\}$ is asymptotically
  compact.
\end{defn}

\begin{lem} \label{lem:Gammalsc}
  If $E_i$ $\Gamma$-converges to $E$, then $E$ is lower
  semi-continuous.
\end{lem}

\begin{proof}
  Let a sequence $x_j \in X$, $j \in \N$, converge to a point $x \in
  X$.  By ($\Gamma$1), for each $j$ there is a net $x_{j,i} \in X_i$
  converging to $x_j$ such that $\lim_i E_i(x_{j,i}) = E(x_j)$.  There is a
  sequence $i(j) \to +\infty$ such that $| E_{i(j)}(x_{j,i(j)}) -
  E(x_j) | < 1/j$ and $x_{j,i(j)} \to x$.  We have
  \[
  \liminf_j E(x_j) = \liminf_j E_{i(j)}(x_{j,i(j)}) \ge E(x).
  \]
  This completes the proof.
\end{proof}

We have the following theorem in the same way as in the linear
version, Lemma 2.14 of \cite{KwSy:specstr}.  The proof is omitted.

\begin{thm} \label{thm:Gamma}
  Assume that $X_i$ and $X$ are all separable.  Then, for any
  $\{E_i\}$ there always exists a $\Gamma$-convergent subnet of
  $\{E_i\}$.
\end{thm}

As an immediate consequence of the theorem, we have:

\begin{cor} \label{cor:asymcpt}
  Assume that $X_i$ and $X$ are all separable.  If $\{E_i\}$ is
  asymptotically compact, it has a compactly convergent subnet.
\end{cor}

The following is a generalization of Remark 2.6 of
\cite{KwSy:specstr}.

\begin{prop} \label{prop:Ecpt}
  Assume that $X$ is separable.  If $E_i$ compactly converges to $E$,
  then $E$ is compact.
\end{prop}

\begin{proof}
  Let $x_j \in X$ be a bounded net with $\limsup_j E(x_j) < +\infty$.
  By ($\Gamma$1), for each $j$ there is a net $x_{j,i} \in X_i$ such
  that, as $i \to \infty$, $x_{j,i} \to x_j$ and $E_i(x_{j,i}) \to
  E(x_j)$.  Take a dense countable subset $\{\xi_n\}_{n \in \N}
  \subset X$ and, for each $n$, a net $\xi_{n,i} \in X_i$ with
  $\xi_{n,i} \to \xi_n$ as $i \to \infty$.  Let us fix a net $N_j \in
  \N$ tending to $\infty$.  By (A3), for any $j$ there is $i(j)$ such
  that $i(j) \to \infty$ as $j \to \infty$,
  \begin{equation}
    \label{eq:cpt}
      |\;d_{X_i}(x_{j,i},\xi_{n,i}) - d_X(x_j,\xi_n)\;| < 1/N_j,
      \quad\text{and}\quad
      |\,E_i(x_{j,i}) - E(x_j)\,| < 1/N_j
  \end{equation}
  for any $i \ge i(j)$ and $n = 1,2,\dots,N_j$.  Since
  \[
  \limsup_j E_{i(j)}(x_{j,i(j)})
  = \limsup_j E(x_j) < +\infty,
  \]
  the asymptotic compactness of $\{E_i\}$ implies that, by replacing
  $\{j\}$ by a sub-directed set, $\{x_{j,i(j)}\}_j$ converges to some
  point $x \in X$.  It follows from \eqref{eq:cpt} that
  \[
  \lim_j d_X(x_j,\xi_n) = d_X(x,\xi_n)
  \]
  for any $n \in \N$.  For any $\epsilon > 0$ there is $n \in \N$ such
  that $d_X(x,\xi_n) < \epsilon$ and hence $d_X(x_j,\xi_n) < \epsilon$
  if $j$ is large enough.  Therefore, by the triangle inequality,
  $d_X(x_j,x) < 2\epsilon$ for all sufficiently large $j$.  This
  completes the proof.
\end{proof}

\begin{defn}[Asymptotic minimizer] \label{defn:min}
  An \emph{asymptotic minimizer of $\{E_i\}$} is defined to be a net
  $x_i \in X_i$ such that $\lim_i (E_i(x_i) - \inf E_i) = 0$.
\end{defn}

An asymptotic minimizer of $\{E_i\}$ always exists whenever $E_i
\not\equiv +\infty$.

\begin{prop} \label{prop:min}
  Assume that $E \not\equiv +\infty$.
  \begin{enumerate}
  \item If $E_i$ $\Gamma$-converges to $E$ and if an asymptotic
    minimizer of $\{E_i\}$ converges to a point $x \in X$, then $x$ is
    a minimizer of $E$.
  \item If $E_i$ compactly converges to $E$, then any bounded
    asymptotic minimizer of $\{E_i\}$ has a subnet converging to a
    minimizer of $E$.
  \end{enumerate}
\end{prop}

\begin{proof}
  We omit the proof of (1) because it is easy.
  
  Let us prove (2).  Since there is a point $y \in X$ with $E(y) <
  +\infty$, ($\Gamma$1) implies that $\inf E_i$ is bounded.  Let $x_i
  \in X_i$ be a bounded asymptotic minimizer of $\{E_i\}$.  Then
  $E_i(x_i)$ is bounded.  By the asymptotic compactness of $\{E_i\}$,
  $\{x_i\}$ has a convergent subnet.  By (1), its limit is a minimizer
  of $E$.  This completes the proof.
\end{proof}

\subsection{Asymptotic compactness under a bound of Poincar\'e constants
  and local covering order}
\label{ssec:Pcov}

Let $M_i$ and $M$ be compact measured metric spaces such that $M_i$
converges to $M$ with respect to the measured Gromov-Hausdorff
topology.  Let $(Y_i,y_i)$ and $(Y,y_0)$ be pointed proper metric
spaces such that $(Y_i,y_i)$ converges to $(Y,y_0)$ in the pointed
Gromov-Hausdorff topology.  By Definition \ref{defn:Lptop}, we have
the $L^p$-topology on $(\bigsqcup_i L^p_{y_i}(M_i,Y_i)) \sqcup
L^p_{y_0}(M,Y)$ for $p \ge 1$.

\begin{rem}
  Since $M_i$ and $M$ are compact, for any maps $\xi_i : M_i \to Y_i$
  and $\xi : M \to Y$ as in Section \ref{ssec:Pcov}, by the
  compactness of $\xi(M)$ we have $L_\xi^p(M,Y) = L_{y_0}^p(M,Y)$ and
  there is $i_0$ such that $L_{\xi_i}^p(M_i,Y_i) = L_{y_i}^p(M_i,Y_i)$
  for any $i \ge i_0$.  The $L^p$-topologies on $(\bigsqcup_{i \ge
    i_0} L^p_{\xi_i}(M_i,Y_i)) \sqcup L^p_{\xi}(M,Y)$ and
  $(\bigsqcup_{i \ge i_0} L^p_{y_i}(M_i,Y_i)) \sqcup L^p_{y_0}(M,Y)$
  coincide.
\end{rem}

Let $E : L^p_{y_0}(M,Y) \to [\,0,+\infty\,]$ be a functional and set
\[
\Dom(E) := \{\,u \in L^p_{y_0}(M,Y) \mid E(u) < +\infty\,\}.
\]
Let $c$, $C$, $\rho$, and $R$ be constants with $c \ge 1$, $C > 0$, and
$0 \le \rho \le R$.
$\mathcal{B}(M)$ denotes the family of Borel subsets of $M$.
Consider the following condition $(P)_{p,c,C,\rho,R}$.

\begin{proclaim}{$(P)_{p,c,C,\rho,R}$:} For any $u \in \Dom(E)$,
  there exists a finitely subadditive function $\mu_u : \mathcal{B}(M)
  \to [\,0,+\infty)$ with $\mu_u(M) \le E(u)$ such that
  \[
  \frac{1}{|B(x,r)|}\iint_{B(x,r)\times B(x,r)} d_Y(u(y),u(z))^p \;dydz
  \le C r^p \, \mu_u(B(x,cr))
  \]
  for any $x \in M$ and $r \in (\,\rho,R\,]$, where
  $|A|$ denotes the measure of $A \subset M$.
\end{proclaim}

The inequality in $(P)_{p,c,C,\rho,R}$ is called a (\emph{weak})
\emph{Poincar\'e inequality}.  Koskela, Shanmugalingam, and Tyson
proved in \cite{KST:Poincare} that $(P)_{2,c,C,0,R}$ for $Y = \R$
implies $(P)_{2,c,C',0,R}$ for general $Y$ under some natural
assumption for $E$.

\begin{defn}[Local covering order] \label{defn:loccov}
  Let $(S,d)$ be a locally compact metric space.  Recall that the
  \emph{order of a covering $\mathcal{U}$ of $S$ at a point $x \in
    S$}, say $\ord_x\mathcal{U}$, is the number of sets in
  $\mathcal{U}$ containing $x$ and that the \emph{order of
    $\mathcal{U}$} is defined to be $\sup_{x \in S}
  \ord_x\mathcal{U}$.  For $c \ge 1$ and $r > 0$, we denote by
  $K_{c,S}(r)$ the maximum of the order of the covering
  $\{B(x_\lambda,cr)\}_\lambda$, where $\{x_\lambda\}$ runs over all
  discrete subsets of $S$ such that $d(x_\lambda,x_{\lambda'}) \ge r$
  for all $\lambda \neq \lambda'$.  For a function $f(r)$, we say that
  \emph{the $c$-local covering order of $S$ is at most $f(r)$} if
  $K_{c,S}(r) \le f(r)$ for any $r > 0$.
\end{defn}

The function $f(r)$ is usually taken as the Landau
symbols $o(r^q)$ and $O(r^q)$, where $\lim_{r \to 0} o(r)/r = 0$ and
$\limsup_{r \to 0} |O(r)|/r < +\infty$.

The proof of the following lemma is straightforward and omitted.

\begin{lem} \label{lem:K}
  If a net of compact metric spaces $S_i$ Gromov-Hausdorff converges
  to a compact metric space $S$, then for any $r > 0$,
  \[
  \limsup_i K_{c,S_i}(r) \le K_{c,S}(r).
  \]
\end{lem}

One of the main results of this paper is the following:

\begin{thm} \label{thm:Pasympcpt}
  Let $p$, $c$, $C$, and $R$ be positive constants with $p,c \ge 1$.
  Assume that the $c$-local covering order of $M$ is at most
  $o(r^{-p})$ as $r \to 0$.  Let $E_i : L^p_{y_i}(M_i,Y_i) \to
  [\,0,+\infty\,]$ be a net of functionals such that each $E_i$
  satisfies $(P)_{p,c,C,\rho_i,R}$ for some net $\rho_i \to 0$ with $0
  \le \rho_i \le R$.  Then, $\{E_i\}$ is asymptotically compact.
\end{thm}

In Section \ref{sec:appl}, we give some interesting applications to
the theorem.

As an immediate consequence of the theorem, we have the following
Rellich compactness result.

\begin{cor} \label{cor:cptemb}
  Let $p$, $c$, $C$, and $R$ be positive constants with $p,c \ge 1$.
  Assume that the $c$-local covering order of $M$ is at most
  $o(r^{-p})$ as $r \to 0$.  Let $E : L^p_{y_0}(M,Y) \to
  [\,0,+\infty\,]$ be a functional satisfying $(P)_{p,c,C,0,R}$.
  Then $E$ is compact.
\end{cor}

Note that the doubling condition implies the boundedness of the
$c$-local covering order, so that our assumption is weaker than the
doubling condition.  Under the doubling condition for the measure on
$M$, Corollary \ref{cor:cptemb} is well-known. According to a result
by Cheeger \cite{Ch:metmeas}, the Poincar\'e inequality and the
doubling condition together imply that $M$ is essentially of finite
dimension.  Theorem \ref{thm:Pasympcpt} and Corollary \ref{cor:cptemb}
hold even for infinite-dimensional spaces, such as,
the space $Q_\infty$ of the following example.

\begin{ex} \label{ex:Q}
  Let us consider the infinite product $Q_\infty := [\,0,1\,] \times
  [\,0,1/2\,] \times \cdots \times [\,0,1/2^n\,] \times \cdots$ with
  $\ell^2$ norm.  This is a compact metric space of infinite
  dimension.  We observe that the $c$-local covering order of
  $Q_\infty$ is at most $O(r^{-\log_2(c+1)})$.  In fact, on $\R^n$ we
  have $K_{c,\R^n}(r) \le (c+1)^n$ by calculating the volumes of
  balls.  Considering $Q_n := [\,0,1\,] \times [\,0,1/2\,] \times
  \cdots \times [\,0,1/2^n\,] \subset \R^n$ we have $K_{c,Q_n}(r) \le
  K_{c,\R^n}(r) \le (c+1)^n$.  There exists a universal constant $a >
  0$ such that $K_{c,Q_\infty}(r) \le a \, K_{c,Q_n}(r)$ if $n \ge
  \log_2(1/r) + a$.  Thus we have $K_{c,Q_\infty}(r) \le
  O(r^{-\log_2(c+1)})$ by taking $n:=[\log_2(1/r)+a+1]$.  Note that if
  $2^p > c+1$ then the assumption for the $c$-local covering order of
  Theorem \ref{thm:Pasympcpt} is satisfied for $M = Q_\infty$.
  
  Let $\mu_n$ be the Lebesgue measure on $Q_n$ normalized as
  $\mu_n(Q_n) = 1$ and $\mu_\infty$ the product probability measure on
  $Q_\infty$.  Then, $Q_n$ measured Gromov-Hausdorff converges to
  $Q_\infty$.  Under the Neumann condition, we define natural
  $p$-energy functionals $\mathcal{E}^{(p)}_n$ and
  $\mathcal{E}^{(p)}_\infty$ on $L^p(Q_n)$ and $L^p(Q_\infty)$
  respectively as follows.  For a smooth function $u : Q_n \to \R$,
  \[
  \mathcal{E}^{(p)}_n(u) :=
  \int_{Q_n} \left( \sum_{k=1}^n
    \left( \frac{\partial u}{\partial x_k} \right)^2
  \right)^{\frac{p}{2}}
  d\mu_n.
  \]
  For a smooth cylindrical function $u : Q_\infty \to \R$ (i.e., $u$
  is a function of the $Q_n$-factor for some $n \in \N$),
  \[
  \mathcal{E}^{(p)}_\infty(u) :=
  \int_{Q_\infty} \left( \sum_{k=1}^\infty
    \left( \frac{\partial u}{\partial x_k} \right)^2
  \right)^{\frac{p}{2}}
  d\mu_\infty,
  \]
  where the sum is in fact a finite sum.  These are both closable (see
  \cite{MR:nonsymDir}) and we denote their closure by the same
  notations $\mathcal{E}^{(p)}_n$ and $\mathcal{E}^{(p)}_\infty$.
  Then, $\mathcal{E}^{(p)}_n$ $\Gamma$-converges to
  $\mathcal{E}^{(p)}_\infty$.  We assume $p = 2$.  Then, the generator
  of $\mathcal{E}^{(2)}_n$ is the Laplacian on $Q_n$.  The $1$-local
  covering order of $M$ is at most $o(r^{-2})$.  We have
  $(P)_{2,1,C,0,R}$ for $\mathcal{E}^{(2)}_n$ for some positive
  constants $C$ and $R$, where $C$ is independent of the dimension $n$
  (see \cite{Bk:isoplog}).  Thus, Theorem \ref{thm:Pasympcpt} implies
  that $\mathcal{E}^{(2)}_n$ compactly converges to
  $\mathcal{E}^{(2)}_\infty$.  This compact convergence is also
  obtained in the following direct way.  The eigenfunctions of the
  Laplacian on $Q_n$ form the products of one-dimensional
  eigenfunctions with Neumann boundary condition.  They produce all
  the eigenfunctions of the generator of $\mathcal{E}^{(2)}_\infty$.
  By this, we can prove that $\mathcal{E}^{(2)}_n$ compactly converges
  to $\mathcal{E}^{(2)}_\infty$.  The compact convergence also holds
  for an infinite Riemannian product $N_1 \times N_2 \times \cdots
  \times N_n \times \cdots$ of closed Riemannian manifolds $N_n$ such
  that the first nonzero eigenvalue of the Laplacian on $N_n$ is
  divergent to infinity as $n \to \infty$.  Since the proof is
  elementary, we omit the details.
\end{ex}

\begin{rem}
  To Theorems 5.1, 5.2, and Corollary 5.1 of \cite{KwSy:specstr}, we
  have to add the assumption that the $c$-local covering order of $X$
  is at most $o(r^{-2})$, which has been missed there.
\end{rem}

The basic strategy of the proof of Theorem \ref{thm:Pasympcpt} is same
as in the linear case in our previous paper \cite{KwSy:specstr}.
However, we need much more delicate discussions and besides recover a
mistake in \cite{KwSy:specstr}.

In \cite{KvS:sobharm}, Korevaar and Schoen proved a Rellich-type
compactness theorem, the idea of which proof seems not to work to
obtain our Theorem \ref{thm:Pasympcpt}.

\begin{proof}[Proof of Theorem \ref{thm:Pasympcpt}]
  We suppose the assumption of the theorem.  As in Section
  \ref{ssec:Lp}, we embed $(Y_i,y_i)$ and $(Y,y_0)$ into a Banach
  space $(\mathbb{B},\|\cdot\|)$ and consider the embeddings
  \eqref{eq:embLp}.  We assume that $y_0 = o$, the origin of
  $\mathbb{B}$.  Let $u_i \in
  L^p_{y_i}(M_i,Y_i)$ be such that
  \[
  \sup_i (E_i(u_i) + d_{L^p}(u_i,y_i)^p) < \infty.
  \]
  Since $y_i \to y_0 = o$ and $|M_i| \to |M| < \infty$, we have
  $\|\iota_i(u_i)-u_i\|_{L^p} = \|y_i\|\,|M_i|^{1/p} \to 0$, where
  $\iota_i$ is the inclusion map in \eqref{eq:embLp}.  Hence,
  $\|u_i\|_{L^p}$ is uniformly bounded.  It suffices to prove that
  $\{u_i\}$ has an $L^p$-convergent subnet in the sense of Definition
  \ref{defn:LptopB}.  Take a sequence of numbers $r_j \searrow 0$, $j
  = 1,2,\dots$ with $r_j \le R$.  For any $i$ and $j$, we find a
  \emph{maximal $r_j$-discrete net} $\{x^i_{jk}\}_{k=1}^{N^i_j}$ of
  $M_i$, i.e., a maximal subset $\{x^i_{jk}\}_{k=1}^{N^i_j}$ of $M_i$
  such that $d_{M_i}(x^i_{jk},x^i_{jk'}) \ge r_j$ for all $k \neq k'$.
  Take a sequence of positive numbers $D_j$ such that
  \begin{equation}
    \label{eq:Dj}
    D_j^p > \sup_{i,k} \frac{2\|u_i\|_{L^p}^p}{|B(x^i_{jk},r_j)|}.
  \end{equation}
  Note that since the measures on $M_i$ and $M$ are of full support,
  $|B(x^i_{jk},r_j)|$ is bounded away from zero for all $i$ and $k$ if
  we fix $j$.  Set 
  \[
  B^i_{jk} := \{\,x \in B(x^i_{jk},r_j) \mid
  \|u_i(x)\| < D_j\,\}
  \]
  and define $\bar u^i_{jk} \in \mathbb{B}$ as the vector-valued
  integral
  \[
  \bar u^i_{jk} := \frac{1}{|B^i_{jk}|} \int_{B^i_{jk}} u_i(x) \; dx.
  \]

  \begin{clm} \label{clm:P}
    For any $i$, $j$, and $k$ with $\rho_i < r_j$, we have
    \[
    \int_{B(x^i_{jk},r_j)} \|u_i(x)-\bar u^i_{jk}\|^p \; dx
    \le 2C r_j^p \,\mu_{u_i}(B(x^i_{jk},cr_j)).
    \]
  \end{clm}

  \begin{proof}
    It follows that
    \[
    \|u_i\|_{L^p}^p \ge \int_{B^i_{jk}} \|u_i(x)\|^p \, dx
    \ge D_j^p \, |B^i_{jk}|,
    \]
    which together with \eqref{eq:Dj} implies $|B^i_{jk}| >
    |B(x^i_{jk},r_j)|/2$.  Therefore, by H\"older's inequality and
    (P), we have
    \begin{align*}
      &\int_{B(x^i_{jk},r_j)} \| u_i(x) - \bar u^i_{jk} \|^p \, dx\\
      &\le \frac{1}{|B^i_{jk}|}
      \iint_{B(x^i_{jk},r_j) \times B^i_{jk}}
      \|u_i(x)-u_i(y)\|^p \, dxdy\\
      &\le \frac{2}{|B(x^i_{jk},r_j)|}
      \iint_{B(x^i_{jk},r_j) \times B(x^i_{jk},r_j)}
      \|u_i(x)-u_i(y)\|^p \, dxdy\\
      &\le 2Cr_j^p \, \mu_{u_i}(B(x^i_{jk},cr_j)).
    \end{align*}
  \end{proof}
  
  Set $U^i_{jk} := B(x^i_{jk},r_j) \setminus \bigcup_{l=1}^{k-1}
  B(x^i_{jl},r_j)$.  Note that $U^i_{jk} \cap U^i_{j l} = \emptyset$
  for any $k \neq l$ and that $\bigcup_{k=1}^{N^i_j} U^i_{jk} = M_i$.
  We define a step map $\bar u^i_j : M_i \to \mathbb{B}$ by $\bar
  u^i_j := \bar u^i_{jk}$ on each $U^i_{jk}$.

  \begin{clm} \label{clm:bar}
    We have
    \[
    \lim_j \limsup_i d_{L^p}(u_i,\bar u^i_j) = 0.
    \]
  \end{clm}

  \begin{proof}
    Claim \ref{clm:P} shows that if $\rho_i < r_j$ then
    \begin{align*}
      d_{L^p}(u_i,\bar u^i_j)^p
      &= \sum_{k=1}^{N^i_j}
      \int_{U^i_{jk}} \|u_i(x)-\bar u^i_{jk}\|^p \; dx
      \le C r_j^p \, \sum_{k=1}^{N^i_j}
      \mu_{u_i}(B(x^i_{jk},cr_j)) \\
      &\le  C r_j^p \, K_{c,M_i}(r_j) \, \mu_{u_i}(M_i).
    \end{align*}
    By Lemma \ref{lem:K} and the assumption, we have
    \[
    \limsup_i K_{c,M_i}(r_j) \le K_{c,M}(r_j) \le o(r_j^{-p}).
    \]
    Since $\mu_{u_i}(M_i) \le E_i(u_i)$ are bounded above, this completes
    the proof of Claim \ref{clm:bar}.
  \end{proof}
  
  For convenience in the later discussions, we define on a metric
  space $(S,d)$ a function $\rho[x,a,b] : S \to [\,0,1\,]$ for $a < b$
  and $x \in S$ by
  \[
  \rho[x,a,b](y) :=
  \begin{cases}
    1 &\text{if $d(x,y) \le a$},\\
    (b-d(x,y))/(b-a) &\text{if $a < d(x,y) < b$},\\
    0 &\text{if $d(x,y) \ge b$},
  \end{cases}
  \]
  for any $y \in S$.  We see that $\rho[x,a,b]$ is a Lipschitz
  function with Lipschitz constant $1/(b-a)$.

  \begin{clm} \label{clm:conv}
    For any fixed $j$ and $k$, $\{\bar u^i_{jk}\}_i$ has a
    convergence subnet in $\mathbb{B}$.
  \end{clm}

  \begin{proof}
    We fix $j$ and $k$.  Let $\epsilon$ be any number with $0 <
    \epsilon < 1$.  Since $u_i(x) \in Y_i \cap B(o,D_j)$ for any $x
    \in B^i_{jk}$, there is a finite sequence of points $p^i_1, p^i_2,
    \dots, p^i_{n(i)} \in Y_i \cap B(o,D_j)$ depending on $\epsilon$
    such that
    \[
    \left\| \bar u^i_{jk} - \frac{1}{n(i)} \sum_{l=1}^{n(i)} p^i_l \right\|
    < \epsilon.
    \]
    There is $i(\epsilon)$ such that $Y_i \cap B(o,D_j) \subset
    B(Y,\epsilon)$ for any $i \ge i(\epsilon)$.  Let $i \ge
    i(\epsilon)$.  We find $q^i_l \in Y$ such that $\|p^i_l - q^i_l\|
    < \epsilon$ and $\|q^i_l\| < D_j + \epsilon$ for any $l$.  By the
    triangle inequalities, we have
    \[
    \left\| \bar u^i_{jk} - \frac{1}{n(i)} \sum_{l=1}^{n(i)} q^i_l \right\|
    < 2\epsilon,
    \]
    so that $\bar u^i_{jk}$ for any $i \ge i(\epsilon)$ is in the
    $2\epsilon$-neighborhood of the closed convex hull, say $C$, of $Y
    \cap \bar B(o,D_j+1)$.  Consequently we have $\lim_i
    d_{\mathbb{B}}(\bar u^i_{jk},C) = 0$.  We apply Mazur's
    compactness theorem which says that the closed convex hull of any
    compact subset of a Banach space is compact
    (see \cite{Mg:introBanach}).  Since $Y \cap \bar B(o,D_j+1)$ is
    compact, $C$ is compact.  This proves Claim \ref{clm:conv}.
  \end{proof}
  
  There is a measure approximation $\{\varphi_i : \Dom(\varphi_i) =
  M_i \to M\}$ such that each $\varphi_i$ is an
  $\epsilon_i$-approximation for some $\epsilon_i \to 0+$.  There is a
  sub-directed set $\mathcal{I}_j$ of $\{i\}$ depending on $j$ such
  that for every $k = 1,\dots,N^i_j$, the limits $x_{jk} := \lim_i
  \varphi_i(x^i_{jk})$, $N_j := \lim_i N^i_j$, and $\bar u_{jk} :=
  \lim_i \bar u^i_{jk}$ all exist.  Replacing with a sub-directed set
  of $\mathcal{I}_j$, we assume that $N_j = N^i_j$ for all $i \in
  \mathcal{I}_j$.  We may also assume that $\mathcal{I}_{j+1} \subset
  \mathcal{I}_j$ for any $j$.  Therefore, by a diagonal argument, we
  find a common cofinal subnet of all $\mathcal{I}_j$ and denote it by
  $\mathcal{I}$.  Set $U_{jk} := B(x_{jk},r_j) \setminus
  \bigcup_{l=1}^{k-1} B(x_{jl},r_j)$.  For any $\epsilon > 0$, $x \in
  M$, and any set $A$, we define
  \begin{align*}
    \chi_{B(x,r_j)}^\epsilon
    &:= \rho[x,r_j-2\epsilon,r_j-\epsilon] : M \to [\,0,1\,], \\
    \chi_{U_{jk}}^\epsilon &:= \chi_{B(x_{jk},r_j)}^\epsilon \cdot
    \prod_{l=1}^{k-1} (1-\chi_{B(x_{jl},r_j)}^\epsilon)
    : M \to [\,0,1\,], \\
    I_A(x) &:=
    \begin{cases}
      1 &\text{if $x \in A$,}\\
      0 &\text{if $x \notin A$.}
    \end{cases}
  \end{align*}

  \begin{clm} \label{clm:chi}
    We have
    \begin{gather}
      \lim_{\epsilon \to 0+}
      \|\chi_{U_{jk}}^\epsilon - I_{U_{jk}}\|_{L^p(M)} = 0,
      \label{eq:chi1} \\
      \lim_{\epsilon \to 0+} \lim_{i \in \mathcal{I}}
      \|\chi_{U_{jk}}^\epsilon \circ \varphi_i -
      I_{U^i_k}\|_{L^p(M_i)} = 0 \label{eq:chi2}
    \end{gather}
    for any $j = 1,2,\dots$ and $k = 1,\dots,N_j$.
  \end{clm}

  \begin{proof}
    Take any $j = 1,2,\dots$ and $k = 1,\dots,N_j$ and fix them.
    Let $\epsilon$ be any positive number.
    We set $A(x,r,r') := B(x,r') \setminus B(x,r)$.  Since $\{I_{U_{jk}}
    \neq \chi_{U_{jk}}^\epsilon\} \subset \bigcup_{l=1}^{N_j}
    A(x_{jl},r_j-2\epsilon,r_j-\epsilon)$, we have
    \[
    \|I_{U_{jk}}-\chi_{U_{jk}}^\epsilon\|_{L^p(M)}^p
    \le \sum_{l=1}^{N_j} |A(x_{jl},r_j-2\epsilon,r_j-\epsilon)|,
    \]
    which implies \eqref{eq:chi1}.
    
    Let $i \in \mathcal{I}$ satisfy $\dis \varphi_i < \epsilon/2$.  For
    any $y \in M_i$ and $l = 1,\dots,N_j$, we have
    \[
    |\,d_{M_i}(x^i_{jl},y)-d_M(x_{jl},\varphi_i(y))\,| < \epsilon/2
    \]
    and hence $\{I_{U^i_{jk}} \neq \chi_{U_{jk}}^\epsilon \circ
    \varphi_i\} \subset \bigcup_{l=1}^{N_j}
    A^i_{l\epsilon}$,
    where $A^i_{l\epsilon} :=
    A(x^i_{jl},r_j-3\epsilon,r_j-\epsilon/2)$, so that
    \[
    \|I_{U^i_{jk}}-\chi_{U_{jk}}^\epsilon\circ
    \varphi_i\|_{L^p(M_i)}^p \le
    \sum_{l=1}^{N_j} |A^i_{l\epsilon}|.
    \]
    Setting
    \[
    \alpha_{l\epsilon} :=
    \rho[x_{jl},r_j-\epsilon/4,r_j-\epsilon/8] \cdot
    (1-\rho[x_{jl},r_j-5\epsilon,r_j-4\epsilon]) : M \to [\,0,1\,],
    \]
    we have $I_{A^i_{l\epsilon}} \le \alpha_{l\epsilon}
    \circ \varphi_i$ for large $i$ and so
    \[
    \limsup_{i \in \mathcal{I}} |A^i_{l\epsilon}| \le
    \limsup_{i \in \mathcal{I}} \int_{M_i} \alpha_{l\epsilon}
    \circ \varphi_i \; dx
    = \int_M \alpha_{l\epsilon} \; dx,
    \]
    which tends to zero as $\epsilon \to 0$.  This completes the proof
    of the claim.
  \end{proof}
  
  For $\bar u_{jk} = \lim_i \bar u^i_{jk}$ we define two maps $\bar
  u_j, \tilde u_j^\epsilon : M \to \mathbb{B}$ by
  \[
  \bar u_j(x) := \sum_{k=1}^{N_j} I_{U_{jk}}(x) \bar u_{jk} 
  \quad \text{and} \quad
  \tilde u_j^\epsilon(x) := \sum_{k=1}^{N_j} \chi_{U_{jk}}^\epsilon(x)
  \, \bar u_{jk}, \quad x \in M.
  \]
  Claim \ref{clm:chi} proves:

  \begin{clm} \label{clm:tilde}
    For any $j = 1,2,\dots$ we have
    \[
      \lim_{\epsilon \to 0+} d_{L^p}(\tilde u_j^\epsilon,\bar
      u_j) = 0
      \quad\text{and}\quad
      \lim_{\epsilon \to 0+} \limsup_i d_{L^p}(\Phi_i \tilde
      u_j^\epsilon,\bar u^i_j) = 0,
    \]
    where $\Phi_i$ is defined by \eqref{eq:Phii}.
    Consequently, $\bar u^i_j$ $L^p$-converges to $\bar u_j$.
  \end{clm}
  
  \begin{clm} \label{clm:Cauchy}
    $\{\bar u_j\}$ is Cauchy in $L^p_o(M,\mathbb{B})$.
  \end{clm}
  
  \begin{proof}
    By Claim \ref{clm:tilde} and (A3), for any $j$ and $j'$,
    \[
    d_{L^p}(\bar u_j,\bar u_{j'})
    = \lim_i d_{L^p}(\bar u^i_j,\bar u^i_{j'})
    \le \limsup_i (d_{L^p}(\bar u^i_j,u_i) 
    + d_{L^p}(u_i,\bar u^i_{j'})),
    \]
    which tends to zero as $j,j' \to \infty$ because of Claim
    \ref{clm:bar}.
  \end{proof}
  
  We set $u := \lim_{j\to\infty} \bar u_j \in L^p_o(M,\mathbb{B})$.
  Since $\bar u^i_j$ $L^p$-converges to $\bar u_j$, by Claim
  \ref{clm:bar}, and by (A4), we obtain that $u_i$ $L^p$-converges to
  $u$.  This completes the proof of Theorem \ref{thm:Pasympcpt}.
\end{proof}

\subsection{Gromov-Hausdorff convergence of energy-sublevel sets}
\label{ssec:GHsublevel}

We give an asymptotic relation between metric spaces $X_i$ and $X$.
Assume in this section that $X_i$ and $X$ are all separable.  Define,
for $c \in \R$,
\[
X^c := \{\,x \in X \mid E(x) \le c\,\}, \qquad
X_i^c := \{\,x \in X_i \mid E_i(x) \le c\,\},
\]
the $c$-sublevel sets of $E$ and $E_i$.
The main purpose of this section is to prove:

\begin{thm} \label{thm:cpt}
  Assume that $E$ and $E_i$ are all compact, i.e., $X^c$ and $X_i^c$
  are all proper for any $c \in \R$.  Then the following (1) and (2)
  are equivalent.
  \begin{enumerate}
  \item $E_i$ compactly converges to $E$.
  \item For any $c \in \R$ there exist a net $c_i \searrow c$ of
    numbers and a net $o_i \in X_i$ of points converging to a point $o
    \in X$ such that the pointed space $(X_i^{c_i},o_i)$ converges to
    $(X^c,o)$ in the Gromov-Hausdorff topology which is compatible
    with the asymptotic relation between $\{X_i\}$ and $X$.
  \end{enumerate}
\end{thm}

The condition (2) of the theorem is corresponding to the spectral
concentration due to Gromov (see Section 3$\tfrac12$.57 of
\cite{Gr:greenbook}).

We give a simple example.

\begin{ex}
  Let $M_i := \{a_i,b_i\}$ be a net of metric spaces consisting of two
  different points $a_i$ and $b_i$ such that $d_{M_i}(a_i,b_i) \to 0$,
  and let $M := \{a\}$ be a space consisting of a single point $a$.
  We equip $M_i$ and $M$ with measures $dx_i := \alpha\,d\delta_{a_i}
  + \beta\,d\delta_{b_i}$ and $dx := (\alpha + \beta)\,d\delta_a$
  respectively, where $\alpha$ and $\beta$ are positive constants and
  $\delta_z$ denotes Dirac's $\delta$-measure at a point $z$.  $M_i$
  converges to $M$ with respect to the measured Gromov-Hausdorff
  topology.  For a pointed proper metric space $(Y,y_0)$, let us
  consider $X_i := L^2_{y_0}(M_i,Y)$ and $X := L^2_{y_0}(M,Y)$.  $X$
  is isometric to $\sqrt{\alpha+\beta}\,Y$, where $c\,Y$ means the
  space $Y$ with metric multiplied by a positive constant $c$.  It
  follows that
  \[
  d_{L^2}(u,v)^2 = \alpha\,d_Y(u(a_i),v(a_i))^2
  + \beta\,d_Y(u(a_i),v(a_i))^2,
  \qquad u, v \in X_i.
  \]
  Therefore, $X_i$ is isometric to the product $\sqrt{\alpha}\,Y
  \times \sqrt{\beta}\,Y$.  We have the map $\varphi_i : M_i \to M$
  defined by $\varphi(a_i) := \varphi(b_i) := a$, which induces a map
  $\Phi_i : X \to X_i$, $\Phi_i(u) := u\circ\varphi_i$.  This
  coincides with the diagonal embedding
  \[
  \sqrt{\alpha+\beta}\,Y
  \hookrightarrow \sqrt{\alpha}\,Y \times \sqrt{\beta}\,Y,
  \quad y \mapsto (y,y).
  \]
  We define natural energy functionals on $X_i$ and $X$ as
  follows:
  \[
  E_i(u) := \frac{d_Y(u(a_i),u(b_i))^2}{d_{M_i}(a_i,b_i)^2},
  \quad\quad E(v) := 0
  \]
  for any $u \in X_i$ and $v \in X$.  Then, we observe that for a
  fixed number $c \ge 0$, $X_i^c$ shrinks to $X^c = X =
  \sqrt{\alpha+\beta} Y$.
\end{ex}

For the proof of the theorem we present some definitions and notations.

\begin{defn}[$\liminf_i A_i$, $\limsup_i A_i$]
  Let $A_i \subset X_i$.  We define two subsets $\liminf_i A_i$ and
  $\limsup_i A_i$ of $X$ by the following (LI) and (LS):
  \begin{itemize}
  \item[(LI)] $x \in \liminf_i A_i$ iff there exists a net $x_i \in A_i$
    converging to $x$.
    \medskip
  \item[(LS)] $x \in \limsup_i A_i$ iff there exist a sub-directed set
    $\{j\}$ of $\{i\}$ and a net $x_j \in A_j$ such that $x_j \to x$.
  \end{itemize}
\end{defn}
It is obvious that $\liminf_i A_i \subset \limsup_i A_i$.

\begin{nota}
  Let $c \in \R$ be any number.
  We define
  \[
  \underline{X}^c := \liminf_i X_i^c \quad\text{and}\quad
  \overline{X}^c := \limsup_i X_i^c.
  \]
  Let $o_i \in X_i$ be a net converging to a point $o \in X$.  For $r
  > 0$, we denote by $\bar B(o,r)$ the closure of the ball $B(o,r)$
  and define
  \begin{alignat*}{2}
    X_i^c(r) &:= X_i^c \cap \bar B(o_i,r), &\qquad
    X^c(r) &:= X^c \cap \bar B(o,r),\\
    \underline{X}^c(r) &:= \underline{X}^c \cap \bar B(o,r), &\qquad
    \overline{X}^c(r) &:= \overline{X}^c \cap \bar B(o,r).\\
  \end{alignat*}
\end{nota}

To prove the implication (1) $\Rightarrow$ (2) of Theorem
\ref{thm:cpt}, we need some lemmas.

\begin{lem} \label{lem:cpt1}
  Under {\rm(1)} of Theorem \ref{thm:cpt}, we have
  \[
  \bigcap_{c' > c} \underline{X}^{c'}
  = \bigcap_{c' > c} \overline{X}^{c'}
  = X^c
  \]
  for any $c \in \R$.
\end{lem}

\begin{proof}
  By using ($\Gamma$1), it is easy to show that
  $\underline{X}^{c'} \supset X^c$ for any $c < c'$, and hence
  \[
  \bigcap_{c' > c} \underline{X}^{c'} \supset X^c.
  \]
  Take any $x \in \bigcap_{c' > c} \overline{X}^{c'}$.
  For the lemma it suffices to show that $x \in X^c$.
  In fact, for any $c' > c$, since $x \in \overline{X}^{c'}$,
  there exist a sub-directed set $\{j\}$ of $\{i\}$ and a net $x_j \in
  X_j^{c'}$ such that $x_j \to x$.  By ($\Gamma$2) we have
  \[
  E(x) \le \liminf_j E_j(x_j) \le c'.
  \]
  By the arbitrariness of $c'$ with $c' > c$, we obtain $x \in X^c$.
  This completes the proof of the lemma.
\end{proof}

\begin{lem} \label{lem:cpt2}
  Under {\rm(1)} of Theorem \ref{thm:cpt}, for any $c \in \R$ and $r > 0$
  $\overline{X}^c(r)$ is compact.
\end{lem}

\begin{proof}
  By Lemma \ref{lem:cpt1}, $\overline{X}^c$ is contained in $X^c$.  It
  follows from a diagonal argument that $\overline{X}^c$ is a closed
  subset of $X$.  Since $X^c$ is proper, $\overline{X}^c(r)$ is
  compact.
\end{proof}

\begin{lem} \label{lem:cpt3}
  Under {\rm(1)} of Theorem \ref{thm:cpt}, for any $c \in \R$, $\epsilon >
  0$ and $r > r_0 > 0$, there exists $c' = c'(c,\epsilon,r,r_0) > c$
  such that
  \[
  X^c(r_0) \subset \underline{X}^{c'}(r_0)
  \subset \overline{X}^{c'}(r)
  \subset B(X^c(r),\epsilon).
  \]
\end{lem}

\begin{proof}
  Lemma \ref{lem:cpt1} implies that $X^c(r_0) \subset
  \underline{X}^{c'}(r_0)$ for any $c' > c$, and
  $\underline{X}^{c'}(r_0) \subset \overline{X}^{c'}(r)$ is obvious.
  Thus, it suffices to prove $\overline{X}^{c'}(r) \subset
  B(X^c(r),\epsilon)$ for some $c'$ with $c' > c$.
  Suppose the contrary, so that there exists a sequence $c_j \searrow
  c$ of numbers and a sequence $x_j \in \overline{X}^{c_j}(r)
  \setminus B(X^c(r),\epsilon)$.  By Lemma \ref{lem:cpt2},
  $x_j$ has a convergent subsequence and we may assume that $x_j$
  converges to a point $x \in X$.  Then $x$ belongs to
  $\bigcap_{j=1}^\infty \overline{X}^{c_j}(r) \setminus
  B(X^c(r),\epsilon)$.  However, by Lemma \ref{lem:cpt1},
  $\bigcap_{j=1}^\infty \overline{X}^{c_j}(r) = X^c$.
  This is a contradiction and completes the proof.
\end{proof}

\begin{lem} \label{lem:cpt4}
  Assume {\rm(1)} of Theorem \ref{thm:cpt} and let $\{f_i : X \to X_i\}$ be
  asymptotically continuous metric approximation compatible with the
  asymptotic relation between $\{X_i\}$ and $X$ (such a $\{f_i\}$
  always exists by Lemma \ref{lem:asympconti}).  Then for any $c' > c$
  and $\epsilon, r > 0$ there exists $i_0 = i_0(c,c',\epsilon,r)$ such
  that
  \[
  X_i^c(r) \subset B(f_i(\overline{X}^{c'}(r)),\epsilon)
  \]
  for any $i \ge i_0$.
\end{lem}

\begin{proof}
  Suppose not.  Then, there exists a net $x_i \in X_i^c(r)$ such that
  \begin{equation}
    \label{eq:cpt1}
    d_{X_i}(x_i,f_i(\overline{X}^{c'}(r))) \ge \epsilon.
  \end{equation}
  Since $\{E_i\}$ is asymptotically compact, $x_i$ has a convergent
  subnet and we denote its limit by $x$.  It follows that $x$ belongs
  to $\bar B(o,r)$.  By recalling the definition of $\overline{X}^c$,
  we have $x \in \overline{X}^c(r)$.  Since $d_{X_i}(x_i,f_i(x)) \to
  0$, this contradicts \eqref{eq:cpt1}.
\end{proof}

\begin{lem} \label{lem:cpt5}
  Under {\rm(1)} of Theorem \ref{thm:cpt}, for any $c' > c$, $\epsilon >0$,
  and $r > r_0 > 0$ there exists $i_1 = i_1(c,c',\epsilon,r,r_0)$ such
  that
  \[
  f_i(\underline{X}^{c'}(r_0)) \subset B(X_i^{c'}(r),\epsilon)
  \]
  for any $i \ge i_1$.
\end{lem}

\begin{proof}
  Suppose not.  Then there exists a point $x \in
  \underline{X}^{c'}(r_0)$ such that $d_{X_i}(f_i(x),X_i^{c'}(r)) \ge
  \epsilon$.  By the definition of $\underline{X}^{c'}$, there exists
  a net $x_i \in X_i^{c'}$ converging to $x$.  Since $d_X(o,x) \le r_0
  < r$, if $i$ is sufficiently large, $x_i$ belongs to $X_i^{c'}(r)$.
  We also have $d_{X_i}(f_i(x),x_i) \to 0$, which is a contradiction.
\end{proof}

\begin{proof}[Proof of Theorem \ref{thm:cpt}]
  Let us first prove (1) $\Rightarrow$ (2). We assume (1).  Let $r_0 >
  0$.  Since $\bigcap_{r > r_0} X^c(r) = X^c(r_0)$, for any $\epsilon
  > 0$ there exists $r > r_0$ such that $d_H(X^c(r),X^c(r_0)) <
  \epsilon$.  Combining this with Lemmas \ref{lem:cpt3},
  \ref{lem:cpt4}, and \ref{lem:cpt5} proves
  \[
  \limsup_i d_{GH}(X^c(r_0),X_i^{c'}(r)) < 4\epsilon.
  \]
  Therefore, by taking a sequence $\epsilon_j \searrow 0$, there exist
  sequences $r_j \searrow r_0$ and $c_j \searrow c$ such that
  \[
  \limsup_i d_{GH}(X^c(r_0),X_i^{c_j}(r_j)) < \epsilon_j,
  \]
  so that we can find a divergent sequence $i(j)$ such that for any $i \ge
  i(j)$,
  \[
  d_{GH}(X^c(r_0),X_i^{c_j}(r_j)) < \epsilon_j.
  \]
  Therefore, there exists a monotone increasing divergent net
  $\{j(i)\}$ such that
  \[
  d_{GH}(X^c(r_0),X_i^{c_{j(i)}}(r_{j(i)})) < \epsilon_{j(i)},
  \]
  which implies (2).

  Let us prove (2) $\Rightarrow$ (1).
  Assume (2).  To prove the asymptotic compactness of $\{E_i\}$,
  we take a net $x_i \in X_i^c(r)$ for some fixed $c \in \R$ and $r >
  0$.  It suffices to prove that $x_i$ has a convergent subnet.
  In fact, since $(X_i^{c_i},o_i)$ converges to $(X^c,o)$ in the
  Gromov-Hausdorff topology for some $c_i \searrow c$ and since $x_i$
  belongs to $X_i^{c_i}$, it has a convergent subnet.

  To prove ($\Gamma$2), we take a net $x_i \in X_i$ converging to a
  point $x \in X$.  We may assume that $\limsup_i E_i(x_i) <
  +\infty$.  Let $c$ be any number with $c > \liminf_i E_i(x_i)$.
  There exists a subnet $\{j\}$ of $\{i\}$ such that
  $\liminf_i E_i(x_i) = \lim_j E_j(x_j)$ and $E_j(x_j) < c$ for any
  $j$.  By the assumption, $(X_j^{c_j},o_i)$ converges to $(X^c,o)$
  for some net $c_j \searrow c$.  Since $x_j$ belongs to $X_j^{c_j}$, $x$
  belongs to $X^c$.  By the arbitrariness of $c$, we have
  $E(x) \le \liminf_i E_i(x_i)$.
  
  Finally we shall prove ($\Gamma$1).  Let $x \in X$.  If $E(x) =
  +\infty$, then (A2) of Definition \ref{defn:asymprel} tells us to
  find a net $x_i \in X_i$ converging to $x$ and we have $E_i(x_i) \to
  E(x)$ by ($\Gamma$2).  Assume $E(x) < +\infty$ and set $c := E(x)$.
  Since there exists a net $c_i \searrow c$ such that
  $(X_i^{c_i},o_i)$ converges to $(X^c,o)$, we can find a net $x_i \in
  X_i^{c_i}$ converging to $x$.  Then we have $E_i(x_i) \le c_i
  \searrow c = E(x)$, which together with ($\Gamma$2) shows $E(x_i)
  \to E(x)$.  This completes the proof of Theorem \ref{thm:cpt}.
\end{proof}

\begin{defn}[Essentially critical value]
  A value $c$ of $E$ is said to be \emph{essentially critical} if the
  closure of the set $\{\,x \in X \mid E(x) < c\,\}$ is a proper subset
  of $X^c$.
\end{defn}

A metric space is called a \emph{geodesic space} if any two points in
the space is joined by a length-minimizing curve, called a
\emph{geodesic}, and the distance between them coincides with the
length of the geodesic.

\begin{prop}
  If $X$ is a geodesic space and if $E : X \to [\,0,+\infty\,]$ is a
  convex function, then any non-minimal value of $E$ is not
  essentially critical.
\end{prop}

\begin{proof}
  Let $c \in \R$ be a non-minimal value of $E$ and take a point $x \in
  X$ with $E(x) = c$.  It suffices to prove that there exists a
  sequence in $\{E < c\}$ converging to $x$.  In fact, since $c$ is
  not minimal, we can find a point $y \in X$ with $E(y) < c$ and join
  $x$ and $y$ by a geodesic segment $\gamma$.  It follows from the
  convexity of $E$ that $\gamma \setminus \{x\}$ is contained in $\{E
  < c\}$.  This completes the proof.
\end{proof}

\begin{thm} \label{thm:cpt2}
  Assume that $E$ and $E_i$ are all compact.  If $E_i$ compactly
  converges to $E$ and if a number $c \in \R$ is not an essentially
  critical value of $E$, then for any net $c_i \to c$ (not necessarily
  monotone), $(X_i^{c_i},o_i)$ converges to $(X^c,o)$ in the
  Gromov-Hausdorff topology compatible with the asymptotic relation
  between $\{X_i\}$ and $X$, where $o_i \in X_i$ is any net converging
  to a point $o \in X$.
\end{thm}

\begin{proof}
  Since $c$ is not essentially critical, we have
  \[
  X^c = \overline{\{E < c\}} = \overline{\bigcup_{\delta > 0}
    X^{c-\delta}},
  \]
  and hence, as $\delta \searrow 0$, $X^{c-\delta}$ converges to $X^c$
  in the compact Hausdorff convergence.
  On the other hand, since
  \[
  X^c = \bigcap_{\delta > 0} X^{c+\delta},
  \]
  as $\delta \searrow 0$, $X^{c+\delta}$ converges to $X^c$ in the
  compact Hausdorff convergence.  Applying Theorem \ref{thm:cpt}
  yields that for any $\delta > 0$ there exist two nets $\rho_i
  \searrow 0$ and $\rho_i' \searrow 0$ such that
  $(X_i^{c-\delta+\rho_i},o_i) \to (X^{c-\delta},o)$ and
  $(X_i^{c+\delta+\rho_i'},o_i) \to (X^{c+\delta},o)$ in the
  Gromov-Hausdorff topology, where $o_i \in X$ is a net converging to
  a point $o \in X$.  Let $c_i$ be any net of numbers converging to
  $c$.  For all sufficiently large $i$, we have $c-\delta+\rho_i < c_i
  < c+\delta+\rho_i'$ and so $X_i^{c-\delta+\rho_i} \subset X_i^{c_i}
  \subset X_i^{c+\delta+\rho_i'}$.  Therefore, $(X_i^{c_i},o_i)$
  converges to $(X^c,o)$ in the Gromov-Hausdorff topology.
\end{proof}

\section{Variational convergence over $\CAT(0)$-spaces}
\label{sec:CAT0}

\subsection{$\CAT(0)$-spaces}

A \emph{geodesic} in a metric space is a rectifiable curve joining two
points $x$ and $y$ in the metric space that is length-minimizing among
all curves joining $x$ and $y$.  Assume that all geodesics have
parameters proportional to arclength.  A \emph{geodesic space} is a
metric space any two points of which can be joined by a
\emph{geodesic} and the distance between them coincides with the
length of the geodesic.  Note that a geodesic joining two fixed points
is not necessarily unique.  We denote by $xy$ one of geodesics joining
$x$ to $y$.  A \emph{triangle $\triangle xyz$} in a geodesic space
consists of three edges $xy$, $yz$, and $zx$.  For any triangle
$\triangle xyz$ there is a triangle $\triangle \tilde x\tilde y\tilde
z$ in $\R^2$ having same side lengths as $\triangle xyz$.  We call
such a $\triangle \tilde x\tilde y\tilde z$ a \emph{comparison
  triangle of $\triangle xyz$}.  A geodesic space $(H,d_H)$ is called
a \emph{$\CAT(0)$-space} if the following \emph{triangle comparison
  condition} is satisfied.
\begin{itemize}
\item[(T)] Take any triangle $\triangle xyz$ in $H$ and a comparison
  triangle $\triangle \tilde x\tilde y\tilde z$ of it.  For any points
  $a \in xy$, $b \in xz$, $\tilde a \in \tilde x\tilde y$, and $\tilde
  b \in \tilde x\tilde z$ with $d_H(x,a) = d_H(\tilde x,\tilde a)$ and
  $d_H(x,b) = d_H(\tilde x,\tilde b)$, we have $d_H(a,b) \le
  d_H(\tilde a,\tilde b)$.
\end{itemize}

Let $H$ be a complete $\CAT(0)$-space.  It follows from (T) that for
given two points $x, y \in H$, a geodesic $xy$ is unique.  However
geodesics in $H$ may branch in general.  Let $C \subset H$ be a closed
convex set.  Then for any $x \in H$ there exists a unique point in $C$
nearest to $x$, which we denote by $\pi_C(x)$.  The map $\pi_C : H \to
C$ is $1$-Lipschitz (see II.2 of \cite{BH:nonpos}).

If $M$ is a measure space and if $Y$ is $\CAT(0)$ then $L^2_\xi(M,Y)$
is $\CAT(0)$ for any measurable map $\xi : M \to Y$ (see
\cite{Jt:nonpos}).

Assume that an asymptotic relation between metric spaces $\{X_i\}$
and $X$ is given.  Consider the following condition:
\begin{itemize}
\item[(G)] If $\gamma_i : [\,0,1\,] \to X_i$ and $\gamma : [\,0,1\,]
  \to X$ are geodesics such that $\gamma_i(0) \to \gamma(0)$ and
  $\gamma_i(1) \to \gamma(1)$, then $\gamma_i(t) \to \gamma(t)$ for
  any $t \in [\,0,1\,]$.
\end{itemize}

\begin{prop} \label{prop:G}
  \begin{enumerate}
  \item If {\rm(G)} is satisfied and if each $X_i$ is geodesic
    {\rm(}resp.~$\CAT(0)${\rm)}, then so is $X$.
  \item If each $X_i$ is $\CAT(0)$ and $X$ is geodesic, then
    {\rm(G)} is satisfied and $X$ is $\CAT(0)$.
  \end{enumerate}
\end{prop}

\begin{proof}
  We omit (1) and prove only (2).  Let $\gamma_i : [\,0,1\,] \to X_i$
  and $\gamma : [\,0,1\,] \to X$ be geodesics such that $\gamma_i(0)
  \to \gamma(0)$ and $\gamma_i(1) \to \gamma(1)$.  By (A2), for a
  given $t \in [\,0,1\,]$ there exists a net $x_i \in X_i$ converging
  to $\gamma(t)$.  It follows from (A3) that $d_{X_i}(\gamma_i(0),x_i)
  \to tl$ and $d_{X_i}(x_i,\gamma_i(1)) \to (1-t)l$, where $l :=
  d_X(\gamma(0),\gamma(1))$.  Applying the triangle comparison
  condition, (T), to $\triangle x_i\gamma_i(0)\gamma_i(1)$ shows that
  $d_{X_i}(x_i,\gamma_i(t)) \to 0$.  Thus, by (A4), we have
  $\gamma_i(t) \to \gamma(t)$.  We have obtained (G).  By using (1),
  $X$ is $\CAT(0)$.
\end{proof}

\subsection{Weak convergence}

Let $H_i$ and $H$ be complete separable $\CAT(0)$-spaces that have an
asymptotic relation.  We have (G) by Proposition \ref{prop:G}(2).  In
what follows, \emph{strong convergence} means convergence with respect
to the asymptotic relation.

\begin{defn}[Weak convergence]
  We say that a net \emph{$x_i \in H_i$ weakly converges to a point $x
    \in H$} if for any net of geodesic segments $\gamma_i$ in $H_i$
  strongly converging to a geodesic segment $\gamma$ in $H$ with
  $\gamma(0) = x$, $\pi_{\gamma_i}(x_i)$ strongly converges to $x$.
\end{defn}
It is easy to prove that a strong convergence implies a weak
convergence and that a weakly convergent net always has a unique weak
limit.

\begin{lem} \label{lem:weak}
  Assume that $H_i \ni x_i \to x \in H$ weakly and
  $H_i \ni y_i \to y \in H$ strongly.  Then, we have
  \begin{enumerate}
  \item $d_H(x,y) \le \liminf d_{H_i}(x_i,y_i)$;
  \item $d_{H_i}(x_i,y_i) \to d_H(x,y)$ iff
    $x_i \to x$ strongly.
  \end{enumerate}
\end{lem}

\begin{proof}
  (1): We set $\gamma := xy$.  Take a net $\hat x_i \in H_i$
  strongly converging to $x$ and set $\gamma_i := \hat x_i y_i$.
  It follows from the assumption that $\pi_{\gamma_i}(x_i) \to
  x$.  We have
  \[
  \liminf d_{H_i}(x_i,y_i) \ge \liminf
  d_{H_i}(\pi_{\gamma_i}(x_i),y_i)
  = d_H(x,y).
  \]
  
  (2): The part of `if' is obvious.  We shall prove the part of `only
  if'.  Since the equalities in the proof of (1) hold, we have
  $d_{H_i}(x_i,\pi_{\gamma_i}(x_i)) \to 0$ by the triangle comparison.
  Hence, the limits of $x_i$ and $\pi_{\gamma_i}(x_i)$ should coincide
  to each other and is $x$.
\end{proof}

\begin{lem} \label{lem:bddweaksub}
  Let $x_i \in H_i$ be a net and
  $\gamma_i,\sigma_i : [\,0,1\,] \to H_i$ geodesic segments such that
  \[
  \lim_i d_{H_i}(\gamma_i(0),\sigma_i(0)) = \lim_i
  d_{H_i}(\gamma_i(1),\sigma_i(1)) = 0.
  \]
  Then we have
  \[
  \lim_i d_{H_i}(\pi_{\gamma_i}(x_i),\pi_{\sigma_i}(x_i)) = 0.
  \]
\end{lem}

\begin{proof}
  It follows from the assumption and the convexity of $[\,0,1\,] \ni t
  \mapsto d_{H_i}(\gamma_i(t),\sigma_i(t))$ (see II.2 of
  \cite{BH:nonpos}) that
  \begin{equation}
    \label{eq:bddweaksub1}
      \lim_i \sup_t d_{H_i}(\gamma_i(t),\sigma_i(t)) = 0.
  \end{equation}
  We set $y_i := \pi_{\gamma_i}(x_i)$, $z_i := \pi_{\sigma_i}(x_i)$,
  and take $s_i, t_i \in [\,0,1\,]$ as to satisfy $\gamma_i(s_i) =
  y_i$ and $\sigma_i(t_i) = z_i$.  \eqref{eq:bddweaksub1} leads to
  \[
  \lim_i d_{H_i}(\sigma_i(s_i),y_i) =
  \lim_i d_{H_i}(\gamma_i(t_i),z_i) = 0.
  \]
  Since $d_{H_i}(x_i,y_i) \le d_{H_i}(x_i,\gamma_i(t_i))$ and
  $d_{H_i}(x_i,z_i) \le d_{H_i}(x_i,\sigma_i(s_i))$, we have
  \begin{align*}
    \lim_i |d_{H_i}(x_i,y_i) - d_{H_i}(x_i,z_i)| &= 0,\\
    \lim_i |d_{H_i}(x_i,y_i) - d_{H_i}(x_i,\gamma_i(t_i))| &= 0.
  \end{align*}
  By the triangle comparison,
  \[
  d_{H_i}(\gamma_i(t_i),y_i)^2
  + d_{H_i}(x_i,y_i)^2
  \le d_{H_i}(x_i,\gamma_i(t_i))^2.
  \]
  Therefore,
  $d_{H_i}(\gamma_i(t_i),y_i) \to 0$ and so
  $d_{H_i}(y_i,z_i) \to 0$.
\end{proof}

\begin{lem} \label{lem:bddweak}
  Any bounded sequence $x_i \in H_i$ has a weakly convergent
  subnet.
\end{lem}

\begin{proof}
  The proof is based on that by Jost \cite{Jt:equilib}.  We may assume
  that $\{x_i\}$ is a countable sequence, i.e., $i = 1,2,\dots$.  Take
  a dense countable subset $\{\xi_\nu\}_{\nu\in\N} \subset H$.  Let
  $y_0 \in H$ and $y_{0,i} \in H_i$ be such that $y_{0,i} \to y_0$
  strongly.  For each $\nu \in \N$, we take a sequence $\xi_{\nu,i}
  \in H_i$ converging to $\xi_\nu$.  Then we have $\gamma_{\nu,i}^0 :=
  y_{0,i}\xi_{\nu,i} \to \gamma_\nu^0 := y_0\xi_\nu$ ($i \to \infty$).
  Since $d_{H_i}(x_i,y_{0,i})$ is bounded, $w_{\nu,i}^0 :=
  \pi_{\gamma_{\nu,i}^0}(x_i)$ has a convergent subsequence whose
  limit, say $w_\nu^0$, is a point in $\gamma_\nu^0$.  By a diagonal
  argument, we can choose a common subsequence of $\{i\}$ independent
  of $\nu$ for which $w_{\nu,i}^0$ converges to $w_\nu^0$ as $i \to
  \infty$ and denote the subsequence by the same notation $\{i\}$.  We
  take a sequence $\epsilon_m \to 0+$, $m = 0,1,2,\dots$.  Let us
  define $y_{m,i}$ and $y_m$ inductively as follows.  Suppose that
  $y_{m,i}$ and $y_m$ are defined and satisfy $\lim_i y_{m,i} = y_m$.
  We shall define $y_{m+1}$ and $y_{m+1,i}$.  By setting
  $\gamma_{\nu,i}^m := y_{m,i}\xi_{\nu,i}$, the sequence $w_{\nu,i}^m
  := \pi_{\gamma_{\nu,i}^m}(x_i)$ has a convergent subsequence, which
  can be chosen to be independent of $\nu$ by a diagonal argument.  We
  replace the sequence $\{i\}$ by such a subsequence and set $w_\nu^m
  := \lim_i w_{\nu,i}^m$.  There exists a number $\nu(m+1) \in \N$
  such that
  \[
  d_H(y_m,w_{\nu(m+1)}^m) > \sup_{\nu\in\N} d_H(y_m,w_\nu^m) -
  \epsilon_m.
  \]
  Define $y_{m+1} := w_{\nu(m+1)}^m$ and $y_{m+1,i} :=
  w_{\nu(m+1),i}^m$.
  This defines sequences $y_{m,i}$ and $y_m$, $m,i = 1,2,\dots$, such
  that $\lim_i y_{m,i} = y_m$ for each $m$.

  By the triangle comparison, we have
  \[
  d_{H_i}(y_{m,i},y_{m+1,i})^2 + d_{H_i}(y_{m+1,i},x_i)^2 \le
  d_{H_i}(y_{m,i},x_i)^2.
  \]
  Taking a subsequence of $\{i\}$ again, we assume that, for each $m$,
  $d_{H_i}(y_{m,i},x_i)$ converges to some number, say $\lambda_m$.
  Since $d_H(y_m,y_{m+1})^2 + \lambda_{m+1}^2 \le \lambda_m^2$,
  $\{\lambda_m^2\}$ is monotone non-increasing and
  $d_H(y_m,y_{m+1}) \to 0$.  For any $\nu \in \N$,
  \begin{equation}
    \label{eq:1}
    \lim_i d_{H_i}(y_{m,i},w_{\nu,i}^m)
    = d_H(y_m,w_\nu^m)
    < \epsilon_m + d_H(y_m,y_{m+1}) =: \epsilon_m' \to 0.
  \end{equation}
  There exists $\{\nu(l,i)\}$ such that $\lim_i \xi_{\nu(l,i)} = y_l$.
  By \eqref{eq:1}, if $\nu \ll i$ then $d_{H_i}(y_{m,i},w_{\nu,i}^m)
  \le \epsilon_m'$.  We may assume that $\nu(l,i) \ll i$.
  Then we have
  \[
  \limsup_i
  d_{H_i}(y_{m,i},w_{\nu(l,i),i}^m) \le \epsilon_m'.
  \]
  Since $\lim_i \gamma_{\nu(l,i),i}^m = y_my_l$,
  we can choose a common subsequence $\{i\}$ independent of $m$ and $l$
  for which
  $w_{\nu(l,i),i}^m \in \gamma_{\nu(l,i),i}^m$ converges to some point
  on $y_my_l$, say $x_{m,l}$.
  By the above inequality,
  \[
  d_H(y_m,x_{m,l}) \le \epsilon_m'.
  \]
  It follows from Lemma \ref{lem:bddweaksub} that
  $\pi_{y_{m,i}y_{l,i}}(x_i) \to x_{m,l}$ as $i \to \infty$.
  Since $\pi_{y_{m,i}y_{l,i}}(x_i) =
  \pi_{y_{l,i}y_{m,i}}(x_i)$ we have $x_{m,l} =
  x_{l,m}$.
  Therefore,
  \[
  d_H(y_l,y_m) \le d_H(y_l,x_{l,m}) + d_H(x_{l,m},y_m)
  \le \epsilon_l' + \epsilon_m'
  \]
  and $\{y_m\}$ is a Cauchy sequence.  Let $x := \lim_m y_m$.
  By \eqref{eq:1}, for any geodesic $\gamma$ emanating from $x$, there
  exists a sequence of geodesic segments $\gamma_i \subset H_i$
  strongly converging to
  $\gamma$ such that
  $\pi_{\gamma_i}(x_i) \to x$ strongly.
  Moreover, it follows from Lemma \ref{lem:bddweaksub} that $x_i \to
  x$ weakly.  This completes the proof.
\end{proof}

\begin{rem}
  A weakly convergent net $x_i \in H_i$ is not necessarily bounded.
  For example, let $\gamma_i$, $i = 1,2,\dots$, be an infinite
  sequence of different rays in $\R^2$ emanating from the origin $o$,
  and let $H$ be the subspace of $\R^2$ consisting of the union of all
  $\gamma_i$.  We assume that $H$ is equipped with its intrinsic
  (geodesic) metric.  Then, $H$ is a $\CAT(0)$-space (actually a
  tree).  The sequence $\gamma_i(i)$ is unbounded, but weakly
  converges to the origin $o$
\end{rem}

\subsection{Mosco convergence and resolvent}

We give functions $E_i : H_i \to [\,0,+\infty\,]$ and $E : H \to
[\,0,+\infty\,]$.

\begin{defn}[Mosco convergence] \label{defn:mosco}
  We say that \emph{$E_i$ converges to $E$ in the Mosco sense} if
  both ($\Gamma$1) in Definition \ref{defn:Gamma} and the following
  ($\Gamma$2') hold.
  \begin{itemize}
  \item[($\Gamma$2')] If $H_i \ni x_i \to x \in H$ weakly, then
    $E(x) \le \liminf E_i(x_i)$.
  \end{itemize}
\end{defn}
Note that ($\Gamma$2') is a stronger condition than ($\Gamma$2), so that a
Mosco convergence implies a $\Gamma$-convergence.

Jost's definition of Mosco convergence in \cite{Jt:nlinDir} seems not
fitting in view of original Mosco's one.  We adopt the original
definition of Mosco convergence in Definition \ref{defn:mosco}.

It is easy to prove the following proposition.
The proof is omitted.

\begin{prop} \label{prop:Mosco}
  Assume that $\{E_i\}$ is asymptotically compact.
  Then the following {\rm(1)--(3)} are all equivalent to each other.
  \begin{enumerate}
  \item $E_i$ converges to $E$ in the Mosco sense.
  \item $E_i$ $\Gamma$-converges to $E$.
  \item $E_i$ compactly converges to $E$.
  \end{enumerate}
\end{prop}

\begin{defn}[Moreau-Yosida approximation and resolvent, \cite{Jt:convex}]
  For $E : H \to [\,0,+\infty\,]$ we define $E^\lambda : H \to
  [\,0,+\infty\,]$ by
  \[
  E^{\lambda}(x) := \inf_{y \in H} (\lambda E(y) + d_H(x,y)^2),
  \qquad x \in H,\ \lambda > 0,
  \]
  and call it the \emph{Moreau-Yosida approximation} of $E$.  If $E$
  is lower semi-continuous and convex and if $E \not\equiv +\infty$,
  then for any $x \in H$ there exists a unique point, say
  $J^E_\lambda(x) \in H$, such that
  \[
  E^{\lambda}(x) = \lambda E(J^E_\lambda(x)) +
  d_H(x,J^E_\lambda(x))^2.
  \]
  This defines a map $J^E_\lambda : H \to H$, called the
  \emph{resolvent of $E$}.
\end{defn}

Note that if $H$ is a Hilbert space and if $E$ is a closed densely
defined quadratic form on $H$, then we have $J^E_\lambda = (I +
\lambda A)^{-1}$.  Here, $I$ is the identity operator and $A$ the
infinitesimal generator associated with $E$, i.e., the self-adjoint
operator on $H$ such that $\Dom(E) = \sqrt{A}$ and $E(x) =
(\sqrt{A}x,\sqrt{A}x)_H$ for any $x \in \Dom(E)$, where
$(\cdot,\cdot)_H$ is the Hilbert inner product on $H$.

\begin{defn}[Strong convergence of maps]
  We say that a \emph{net of maps $f_i : H_i \to H_i$ strongly
    converges to a map $f : H \to H$} if $f_i(x_i)$ strongly converges
  to $f(x)$ for any strongly convergent net $H_i \ni x_i \to x \in H$.
\end{defn}

\begin{lem} \label{lem:convex}
  Let $E_i : H_i \to [\,0,+\infty\,]$ be lower semi-continuous and
  convex.  If $E_i$ $\Gamma$-converges to $E$, then $E$ is lower
  semi-continuous and convex.
\end{lem}

\begin{proof}
  The lower semi-continuity of $E$ follows from Lemma \ref{lem:Gammalsc}.
  We prove the convexity of $E$.
  Let $\gamma : [\,0,1\,] \to X$ be a given geodesic.  There are two
  nets $x_i,y_i \in X_i$ converging to $\gamma(0), \gamma(1)$
  respectively such that $E_i(x_i) \to E(\gamma(0))$ and $E_i(y_i) \to
  E(\gamma(1))$.  For each $i$ we take a geodesic $\gamma_i :
  [\,0,1\,] \to X_i$ joining $x_i$ to $y_i$.  Let $t \in [\,0,1\,]$ be
  any number.  The convexity of $E_i$ says that
  \[
  E_i(\gamma_i(t))) \le
  (1-t)E_i(\gamma_i(0)) + t E_i(\gamma_i(1)).
  \]
  Since $\gamma_i(t)$ converges to $\gamma(t)$ and by ($\Gamma$2), we
  have
  \[
  E(\gamma(t)) \le \liminf_i E_i(\gamma_i(t)))
  \le (1-t)E(\gamma(0)) + t E(\gamma(1)).
  \]
  This completes the proof.
\end{proof}

In what follows, we assume that functions $E_i : H_i \to
[\,0,+\infty\,]$ and $E : H \to [\,0,+\infty\,]$ are all \emph{lower
  semi-continuous, convex, and are not identically equal to
  $+\infty$}.  Let $J^i_\lambda$ and $J_\lambda$ be the resolvents of
$E_i$ and $E$ respectively.

\begin{prop} \label{prop:resolvent1}
  If $E_i$ converges to $E$ in the Mosco sense, then for any $\lambda
  > 0$ we have the following {\rm(1)} and {\rm(2)}.
  \begin{enumerate}
  \item $E_i^\lambda$ strongly converges to $E^\lambda$.
  \item $J^i_\lambda$ strongly converges to $J_\lambda$.
  \end{enumerate}
\end{prop}

\begin{proof}
  Assume that $E_i$ converges to $E$ in the Mosco sense and let
  $H_i \ni x_i \to x \in H$.  We set $y_i := J^i_\lambda(x_i)$ and $y
  := J_\lambda(x)$.  Since $E \not\equiv +\infty$, we have
  $E^{\lambda}(x) < +\infty$ and hence $E(y) < +\infty$.
  We shall prove that
  \begin{equation}
    \label{eq:resolvent}
    \limsup_i E^{\lambda}(x_i) \le E^\lambda(x)
    = \lambda E(y) + d_H(x,y)^2 < +\infty.
  \end{equation}
  In fact, by ($\Gamma$2), there exists a net $z_i \in H_i$ such that
  $z_i \to y$ and $E_i(z_i) \to E(y)$.
  It follows that $E^{\lambda}(x_i) \le \lambda E(z_i) +
  d_{H_i}(x_i,z_i)^2$, the right-hand side of which converges to
  $\lambda E(y) + d_H(x,y)^2$.  This proves \eqref{eq:resolvent}.
  
  Since $E^{\lambda}(x_i) = \lambda E(y_i) + d_{H_i}(x_i,y_i)^2$, this
  and \eqref{eq:resolvent} together imply that $\{y_i\}$ is bounded,
  so that it has a weakly convergent subnet.  We replace $\{y_i\}$
  with the weakly convergent subnet and denote its weak limit by $y'$.
  By ($\Gamma$2') and Lemma \ref{lem:weak} we have
  \begin{align*}
    \lambda E(y) + d_H(x,y)^2 &\le \lambda E(y') + d_H(x,y')^2 \\
    &\le \liminf_i (\lambda E(y_i) + d_H(x_i,y_i)^2).
  \end{align*}
  By \eqref{eq:resolvent}, the inequalities above become equalities
  and so we have $y = y'$ and $d_{H_i}(x_i,y_i) \to d_H(x,y)$.  Thus,
  $y_i$ strongly converges to $y$ and $E^\lambda(x_i) \to E^\lambda(x)$.
\end{proof}

\begin{prop} \label{prop:resolvent2}
  If $E_i^\lambda$ strongly converges to $E^\lambda$ for any $\lambda
  > 0$, then $E_i$ $\Gamma$-converges to $E$.
\end{prop}

\begin{proof}
  The proof is essentially same as in Lemma 1.4.3 of \cite{Jt:nlinDir}.
  Assume that $E_i^\lambda$ strongly converges to $E^\lambda$ for any
  $\lambda > 0$.  Let a net $x_i \in H_i$ converge to a point $x \in
  H$.  Then we have
  \[
  \liminf_i E_i(x_i) \ge \lim_i \frac{1}{\lambda} E_i^\lambda(x_i)
  = \frac{1}{\lambda} E^\lambda(x) \to E(x)
  \quad\text{as $\lambda \to 0+$}.
  \]
  The rest is to prove the existence of a net $y_i \in H_i$ such that
  $y_i \to x$ and $\limsup_i E_i(y_i) \le E(x)$.  We may assume that
  $E(x) < +\infty$.
  It follows from the assumption that for any $\lambda > 0$,
  \[
  E(x) \ge \frac{1}{\lambda} E^\lambda(x)
  = \lim_i \frac{1}{\lambda} E_i^\lambda(x_i).
  \]
  A diagonal argument shows that there exists a net $\lambda_i \to 0+$
  such that if we set $y_i := J^i_{\lambda_i}(x_i)$ then
  \[
  E(x) \ge \lim_i \frac{1}{\lambda_i} E_i^{\lambda_i}(x_i)
  = \lim_i (E_i(y_i) + \frac{1}{\lambda_i} d_{H_i}(x_i,y_i)^2).
  \]
  This implies that $E(x) \ge \limsup_i E_i(y_i)$ and
  $d_{H_i}(x_i,y_i)^2 \to 0$.  Thus, $y_i$ converges to $x$.  This
  completes the proof.
\end{proof}

\begin{rem}
  For quadratic forms in Hilbert spaces, the strong convergence
  $E^\lambda_i \to E^\lambda$ is equivalent to the Mosco convergence
  $E_i \to E$ (see Theorem \ref{thm:convspec}).  We do not know if
  this still holds in the general case.
\end{rem}

\begin{rem}
  The converse of Proposition \ref{prop:resolvent2} does not hold.  In
  fact, let $H$ be a Hilbert space and $\{e_i\}_{i \in \N}$ a complete
  orthonormal basis on $H$.  We define $E_i(x) := \|x-e_i\|_H^2$ and
  $E(x) := \|x\|_H^2 +1$ for any $x \in H$, where $\|\cdot\|_H$
  denotes the Hilbert norm.  Then, it is easy to see that $E_i$
  strongly converges to $E$ and in particular, $E_i$
  $\Gamma$-converges to $E$.  However, we have, for any $\lambda > 0$,
  \begin{alignat*}{2}
    E_i^\lambda(o) &= \frac{\lambda}{1+\lambda},
    &\quad
    J_i^\lambda(o) &= \frac{\lambda}{1+\lambda} e_i,\\
    E^\lambda(o) &= 1,
    &\quad
    J_\lambda(o) &= o.
  \end{alignat*}
  Thus, $E_i^\lambda$ (resp.~$J_i^\lambda$) does not strongly converge
  to $E^\lambda$ (resp.~$J_\lambda$).
\end{rem}

Propositions \ref{prop:Mosco}, \ref{prop:resolvent1}, and
\ref{prop:resolvent2} together imply the following

\begin{cor} \label{cor:resolvent}
  Assume that $\{E_i\}$ is asymptotically compact.  Then, the
  following {\rm(1)} and {\rm(2)} are equivalent.
  \begin{enumerate}
  \item $E_i$ compactly converges to $E$.
  \item $E_i^\lambda$ strongly converges to $E^\lambda$ for any
    $\lambda > 0$.
  \end{enumerate}
\end{cor}

\begin{defn}[Asymptotic compactness of maps]
  We say that \emph{a net of maps $\{f_i : H_i \to H_i\}$ is
    asymptotically compact} if for any bounded net $x_i \in H_i$,
  $\{f_i(x_i)\}$ has a strongly convergent subnet.  We say that
  \emph{$f_i : H_i \to H_i$ compactly converges to a map $f : H
    \to H$} if $\{f_i\}$ is asymptotically compact and strongly
  converges to $f$.
\end{defn}

\begin{prop} \label{prop:EJasympcpt}
  Assume that there exists a bounded net $z_i \in H_i$ such that
  $\{E(z_i)\}$ is bounded.  Then the following {\rm(1)} and {\rm(2)}
  are equivalent.
  \begin{enumerate}
  \item $\{E_i\}$ is asymptotically compact.
  \item $\{J^i_\lambda\}$ is asymptotically compact for any $\lambda
    > 0$.
  \end{enumerate}
\end{prop}

\begin{proof}
  (1) $\implies$ (2): Assume (1).  Let $x_i \in H_i$ be a bounded net
  and $\lambda > 0$ be fixed.  Then, $E_i(J^i_\lambda(x_i)) +
  d_{H_i}(J^i_\lambda(x_i),x_i)^2 \le E_i(z_i) + d_{H_i}(x_i,z_i)^2$
  is bounded.  Hence, the asymptotic compactness of $\{E_i\}$ yields
  that $\{J^i_\lambda(x_i)\}$ has a strongly convergent subnet.  Thus
  we obtain (2).
  
  (2) $\implies$ (1): Assume (2).  Let $x_i \in H_i$ be a bounded net
  such that $\sup_i E_i(x_i) =: C < +\infty$.  To prove (1), it suffices
  to show the existence of a convergent subnet of $\{x_i\}$.  Since
  $\lambda E(J^i_\lambda(x_i)) + d_{H_i}(J^i_\lambda(x_i),x_i)^2 \le
  \lambda E(x_i)$, we have
  \begin{equation} \label{eq:EJasympcpt}
    d_{H_i}(J^i_\lambda(x_i),x_i)^2 \le \lambda C.
  \end{equation}
  By (2), we have a convergent subnet of $J^i_\lambda(x_i)$ depending
  on each $\lambda > 0$.  Take a sequence $\lambda_k \to 0+$.  By a
  diagonal argument, there exists a common subnet of $\{x_i\}$ for
  which $y_{k,i} := J^i_{\lambda_k}(x_i)$ converges for all $k$.  Set
  $y_k := \lim_i y_{k,i}$.  It follows from \eqref{eq:EJasympcpt} that
  $d_{H_i}(y_{k,i},y_{l,i}) \le
  (\sqrt{\lambda_k}+\sqrt{\lambda_l})\sqrt{C}$ and hence $d_H(y_k,y_l)
  \le (\sqrt{\lambda_k}+\sqrt{\lambda_l})\sqrt{C}$, so that $\{y_k\}$
  is a Cauchy sequence in $H$.  Denote the limit of $\{y_k\}$ by $y$.
  Again by a diagonal argument, we can choose a subnet $\{i(k)\}$ of
  $\{i\}$ in such a way that $y_{k,i(k)}$ converges to $y$.  Since
  $d_{H_{i(k)}}(y_{k,i(k)},x_{i(k)}) \le \sqrt{\lambda_k C} \to 0$, the
  sequence $\{x_{i(k)}\}$ converges to $y$.  This completes the proof.
\end{proof}

\begin{rem}
  The existence of $\{z_i\}$ in Proposition \ref{prop:EJasympcpt} is
  necessary.  In fact, if such a net $\{z_i\}$ does not exist, then
  (1) is always true, but (2) does not necessarily hold.
\end{rem}

\begin{cor}
  The following {\rm(1)} and {\rm(2)} are equivalent.
  \begin{enumerate}
  \item $E$ is compact.
  \item $J_\lambda$ is compact for any $\lambda > 0$, i.e.,
    any bounded subset of $H$ is mapped by $J_\lambda$ to
    a relatively compact set.
  \end{enumerate}
\end{cor}

\begin{proof}
  Set $E_i = E$ in Proposition \ref{prop:EJasympcpt}.
\end{proof}

\begin{cor} \label{cor:EJasympcpt}
  The following {\rm(1)} and {\rm(2)} are equivalent.
  \begin{enumerate}
  \item $E_i$ compactly converges to $E$.
  \item $E_i^\lambda$ strongly converges to $E^\lambda$ and
    $\{J^i_\lambda\}$ is asymptotically compact for any $\lambda >
    0$.
  \end{enumerate}  
\end{cor}

\begin{proof}
  (1) $\implies$ (2): We assume (1).  Let us see the existence of a
  bounded net $z_i \in H_i$ such that $\{E_i(z_i)\}$ is bounded.  In
  fact, this follows from the $\Gamma$-convergence of $E_i$ to $E$ and
  $E \not\equiv +\infty$.  Thus, Corollary \ref{cor:resolvent} and
  Proposition \ref{prop:EJasympcpt} imply (2).
  
  (2) $\implies$ (1): Assume (2).  By Corollary \ref{cor:resolvent}
  and Proposition \ref{prop:EJasympcpt}, it suffices to prove the
  existence of a bounded net $z_i \in H_i$ such that $\{E_i(z_i)\}$ is
  bounded.  To see this, we fix a number $\lambda > 0$.  Since $E
  \not\equiv +\infty$, we have a point $x \in H$ with $E^\lambda(x) <
  +\infty$.  Find a net $x_i \in H$ converging to $x$. Then, (2)
  implies that $E^\lambda(x_i) \to E^\lambda(x) < +\infty$.  Setting
  $z_i := J^i_\lambda(x_i)$, we have $E^\lambda(x_i) = \lambda
  E_i(z_i) + d_{H_i}(x_i,z_i)^2$.  Therefore, $\{z_i\}$ and
  $\{E_i(z_i)\}$ are both bounded.  This completes the proof.
\end{proof}

\begin{thm} \label{thm:cptconv}
  The following {\rm(1)} and {\rm(2)} are equivalent.
  \begin{enumerate}
  \item $E_i$ compactly converges to $E + c$ for some constant $c \in
    \R$.
  \item $J^i_\lambda$ compactly converges to $J_\lambda$ for any
    $\lambda > 0$.
  \end{enumerate}
\end{thm}

To prove Theorem \ref{thm:cptconv}, we need the concept of semigroup,
which is studied by Jost and Mayer.

\begin{defn}[Semigroup, \cite{My:gradflow, Jt:nlinDir}]
  Let $E : H \to [\,0,+\infty\,]$ be a lower semi-continuous and
  convex function and $\Dom(E) := \{\,x \in H \mid E(x) < +
  \infty\,\}$.  For any $x \in \overline{\Dom(E)}$ and $t > 0$, there
  exists the limit
  \[
  T^E_t(x) := \lim_{n \to \infty} (J^E_{t/n})^n(x)
  \]
  (see \cite{My:gradflow, Jt:nlinDir}).  Define $T^E_0$ to be the
  identity map.  The family of the maps $T^E_t : H \to H$, $t \ge 0$,
  is called the \emph{semigroup associated with $E$}.
\end{defn}

Define
\[
|\nabla(-E)|(x) := \limsup_{\substack{y \to x\\ y \neq x}}
\frac{-E(y)+E(x)}{d_H(x,y)}, \qquad x \in \Dom(E).
\]

\begin{thm}[Jost-Mayer \cite{My:gradflow, Jt:nlinDir}] \label{thm:JM}
  \begin{enumerate}
  \item For each $t \ge 0$, $T^E_t : \overline{\Dom(E)} \to H$ is
    Lipschitz continuous with Lipschitz constant $1$.  For each $x \in
    \overline{\Dom(E)}$, $t \mapsto T^E_t(x)$ is continuous on
    $[\,0,+\infty\,)$ and locally Lipschitz continuous on
    $(\,0,+\infty\,)$.
  \item For any $s,t \ge 0$ we have
    \[
    T^E_{s+t} = T^E_s \circ T^E_t.
    \]
  \item For any fixed $x \in \overline{\Dom(E)}$, setting $c(t) :=
    T^E_t(x)$, $t \ge 0$, we have
    \begin{align*}
      \lim_{h \to 0+} \frac{d_H(c(t+h),c(t))}{h}
      &= \lim_{h \to 0+}
      \frac{-E(c(t+h))+E(c(t))}{d_H(c(t+h),c(t))}\\
      &= |\nabla(-E)|(c(t)).
    \end{align*}
  \end{enumerate}
\end{thm}

\begin{lem} \label{lem:cptconv}
  Let $E, F : H \to [\,0,+\infty\,]$ be two lower semi-continuous
  convex functions with $E, F \not\equiv +\infty$.  Then the following
  {\rm(1)--(3)} are equivalent.
  \begin{enumerate}
  \item $\Dom(E) = \Dom(F)$ and $E - F$ is constant on $\Dom(E)$.
  \item $J^E_\lambda = J^F_\lambda$ for any $\lambda > 0$.
  \item $\overline{\Dom(E)} = \overline{\Dom(F)}$ and $T^E_t = T^F_t$
    for any $t \ge 0$.
  \end{enumerate}
\end{lem}

\begin{proof}
  The implications (1) $\implies$ (2) $\implies$ (3) are obvious by
  the definitions of the resolvent and semigroup.
  
  Let us prove (3) $\implies$ (1).  Applying Theorem \ref{thm:JM}(3)
  yields that for any $x \in \overline{\Dom(E)}$,
  \[
  E(x) - \inf E
  = \int_0^{+\infty} \{ |\nabla(-E)|(T^E_t(x)) \}^2 \; dt.
  \]
  Theorem \ref{thm:JM}(3) implies that
  $|\nabla(-E)|$ is determined only by the semigroup $T^E_t$.
  Since $T^E_t = T^F_t$ for any $t \ge 0$,
  we have $|\nabla(-E)| = |\nabla(-F)|$, so that
  $E(x) - \inf E = F(x) - \inf F$.
  This completes the proof.
\end{proof}

\begin{proof}[Proof of Theorem \ref{thm:cptconv}]
  (1) $\implies$ (2) follows from Proposition \ref{prop:resolvent1}
  and Corollary \ref{cor:EJasympcpt}.
  
  Let us prove (2) $\implies$ (1).  Assume (2).  Since
  $\{J^i_\lambda\}_i$ is asymptotically compact for any $\lambda > 0$,
  Proposition \ref{prop:EJasympcpt} implies that $\{E_i\}$ is
  asymptotically compact and has a compactly convergent subnet.  It
  suffices to show that the limit of any compactly convergent subnet
  of $\{E_i\}$ coincides with $E$.  Take a compactly convergent subnet
  of $\{E_i\}$ and denote it by the same notation $\{E_i\}$.  Let $F :
  H \to [\,0,+\infty\,]$ be its limit.  By Proposition \ref{prop:Ecpt}
  and Lemma \ref{lem:convex}, $F$ is compact and convex.  By
  Proposition \ref{prop:resolvent1}, $J^i_\lambda$ strongly converges
  to $J^F_\lambda$ for any $\lambda > 0$.  By (2) we have $J^E_\lambda
  = J^F_\lambda$ for any $\lambda > 0$.  Applying Lemma
  \ref{lem:cptconv} yields that $\Dom(E) = \Dom(F)$ and $E - F$ is
  constant.  This completes the proof.
\end{proof}

\begin{prob}
  What can we say about the semigroup $T^{E_i}_t$ of $E_i$
  associated with the convergence of $E_i$?
  See Theorem \ref{thm:convspec} below for the Hilbert case.
\end{prob}

\subsection{Case of Hilbert spaces} \label{ssec:Hil}

Throughout this section, we assume that $H$ and $H_i$ are real Hilbert
spaces that have an asymptotic relation.  Let $A$ and $A_i$ be
selfadjoint operators on $H$ and $H_i$ with their spectral measures
$\mu$ and $\mu_i$ respectively.  Denote by $\{T_t\}_{t \ge 0}$ and
$\{T_t^i\}_{t \ge 0}$ the strongly continuous contraction semigroups
($T_t := e^{-tA}$, $T_t^i := e^{-tA_i}$, $t \ge 0$), and by
$\{J_\lambda\}_{\lambda > 0}$ and $\{J^i_\lambda\}_{\lambda > 0}$ the
strongly continuous resolvents ($J_\lambda := (I + \lambda A)^{-1}$
and $J^i_\lambda := (I + \lambda A_i)^{-1}$, $\lambda > 0$).  We have
a densely defined closed quadratic form $\calE$ on $H$ defined by
$\calE(u,v) := (\sqrt{A}u,\sqrt{A}v)_H$, $u,v \in \Dom(\calE) :=
\Dom(\sqrt{A})$.  We also have $\calE_i$ in the same manner.  We say
that a \emph{continuous function $\varphi : \R \to \R$ vanishes at
  infinity} if $\lim_{|x| \to \infty} f(x) = 0$.

By Lemma \ref{lem:asymprel}(1), we have a compatible, linear metric
approximation $\{f_i\}$ for $\{H_i\}$ and $H$ such that $\Dom(f_i)$ is
the subspace consisting of finite linear combinations of a basis of
$H$.  Therefore, Assumption 2.1 of \cite{KwSy:specstr} is satisfied.
The strong (weak) topology of Definition 2.4 (2.5) of
\cite{KwSy:specstr} is compatible with that of this paper.  Thus, by
the results of Section 2 of \cite{KwSy:specstr}, we obtain the
following:

\begin{thm} \label{thm:convspec}
  The following are all equivalent{\rm :}
  \begin{enumerate}
  \item\label{item:form} $\calE_i \to \calE$ with respect to the
    Mosco topology {\rm(}resp.~$\calE_i \to \calE$
    compactly{\rm)}.
  \item\label{item:resolv} $J^i_\lambda \to J_\lambda$ strongly
    {\rm(}resp.~compactly{\rm)} for some $\lambda > 0$.
  \item\label{item:semigp} $T_t^i \to T_t$ strongly
    {\rm(}resp.~compactly{\rm)} for some $t > 0$.
  \item\label{item:cptsupp} $\varphi(A_i) \to \varphi(A)$
    strongly {\rm(}resp.~compactly{\rm)} for any continuous function
    $\varphi : [\,0,\infty\,) \to \R$ with compact support.
  \item\label{item:Cinfty} $\varphi_i(A_i) \to \varphi(A)$
    strongly {\rm(}resp.~compactly{\rm)} for any net $\{\varphi_i
    : [\,0,\infty\,) \to \R\}$ of continuous functions vanishing at
    infinity which uniformly converges to a continuous function
    $\varphi : [\,0,\infty\,) \to \R$ vanishing at infinity.
  \item\label{item:interv} $\mu_i((\,a,b\,]) \to \mu((\,a,b\,])$
    strongly {\rm(}resp.~compactly{\rm)} for any two real numbers $a <
    b$ which are not in the point spectrum of $A$.
  \item\label{item:specmeas} $(\mu_i u_i,v_i)_{H_i} \to (\mu u,v)_H$
    vaguely for any nets $u_i,v_i \in H_i$ and any $u,v \in H$ such
    that $u_i \to u$ strongly and $v_i \to v$ weakly {\rm(}resp.~$u_i
    \to u$ weakly and $v_i \to v$ weakly{\rm)}.
  \end{enumerate}
\end{thm}

\begin{defn}
  The \emph{strong graph limit $\Gamma_\infty$ of
    $\{A_i\}$} is defined to be the set of pairs
  $(u,v) \in H \times H$ such that there exists a net of vectors
  $u_i \in \Dom(A_i)$ with $u_i \to u$ and $A_i
  u_i \to v$ strongly.
\end{defn}

\begin{thm} \label{thm:convspec2}
  The following are all equivalent{\rm :}
  \begin{enumerate}
  \item $\calE_i \to \calE$ with respect to the Mosco topology.
  \item $\varphi_i(A_i) \to \varphi(A)$ strongly for any net
    $\{\varphi_i : [\,0,\infty\,) \to \R\}$ of bounded continuous
    functions uniformly converging to a bounded continuous function
    $\varphi : [\,0,\infty\,) \to \R$.
  \item The strong graph limit $\Gamma_\infty$ of
    $\{A_i\}$ coincides with the graph of $A$.
  \end{enumerate}
\end{thm}

\begin{prop} \label{prop:semicontspec}
  If $\calE_i \to \calE$ with respect to the Mosco topology, then
  \[
  \sigma(A) \subset \lim_i \sigma(A_i),
  \]
  i.e., for any $\lambda \in \sigma(A)$ there exist $\lambda_i
  \in \sigma(A_i)$ tending to $\lambda$.
\end{prop}

Define $n(I) := \dim\mu(I)H$ and $n_i(I) := \dim\mu_i(I)H$ for a Borel
subset $I \subset \R$.  Note that if $A$ has only discrete spectrum,
then $n(I)$ coincides with the number of the eigenvalues in $I$ of $A$
with multiplicities.

\begin{prop} \label{prop:multi}
  Let $a < b$ be two numbers which are not in the point spectrum of
  $A$.  If $\calE_i \to \calE$ with respect to the Mosco topology, we
  have
  \begin{equation}
    \label{eq:multi1}
    \liminf_i n_i((\,a,b\,]) \ge n((\,a,b\,]).
  \end{equation}
\end{prop}

\begin{thm} \label{thm:cptspectra}
  Assume that $\calE_i \to \calE$ compactly.  Then, for any
  $a,b \in \R \setminus \sigma(A)$ with $a < b$, we have
  $n_i((\,a,b\,]) = n((\,a,b\,])$ for sufficiently large
  $i$.  In particular, the limit set of $\sigma(A_i)$
  coincides with $\sigma(A)$.
\end{thm}

\begin{rem}
  Assume that $\calE_i \to \calE$ compactly.  Then, $\calE$ is compact
  and so $A$ has only discrete spectrum.  Even if $\calE_i$ are not
  necessarily compact, Theorem \ref{thm:cptspectra} shows that the
  bottom of the essential spectrum of $A_i$ is divergent to $+\infty$.
  Thus, for each $k \in \N$, the $k^{th}$ eigenvalue of $A_i$ is
  well-defined if $i$ is large enough compared with $k$.  We shall see
  an example of asymptotically compact $\{E_i\}$ with noncompact $E_i$
  in Section \ref{ssec:approx}.
\end{rem}

\begin{cor} \label{cor:spectra}
  Assume that $\calE_i \to \calE$ compactly.  Denote by $\lambda_k$
  {\rm(}resp.~$\lambda_k^i${\rm)} the $k^{th}$ eigenvalue of $A$
  {\rm(}resp.~$A_i${\rm)} with multiplicity.  {\rm(}$\lambda_k^i$ is
  defined if $i$ is large enough compared with $k$.{\rm)} We
  set $\lambda_k := \infty$ for all $k \ge \dim H + 1$ if $\dim H <
  \infty$, and $\lambda_k^i := \infty$ for all $k \ge \dim H_i + 1$ if
  $\dim H_i < \infty$.  Then we have
  \[
  \lim_i \lambda_k^i = \lambda_k
  \qquad \text{for any $k$.}
  \]
  Moreover, let $\{\varphi_k^i\}_{k=1,2,\dots}$ be a {\rm(}possibly
  incomplete{\rm)} orthonormal basis on $H_i$ such that $\varphi_k^i$
  is an eigenvector for $\lambda_k^i$ of $A_i$.  Then, by replacing
  with a sub-directed set of $\{i\}$ if necessarily, for each fixed $k
  \in \N$ with $k \le \dim H$, the vector $\varphi_k^i$ strongly
  converges to some eigenvector $\varphi_k$ for $\lambda_k$ of $A$
  such that $\{\varphi_k\}_{k=1}^{\dim H}$ is a complete orthonormal
  basis on $H$.
\end{cor}

\section{Application} \label{sec:appl}

\subsection{Compact convergence of approximating energy functional}
\label{ssec:approx}

Let $(M,d_M)$ be a compact measured metric space and $(Y,d_Y)$ a proper
metric space.  Let $p \ge 1$ be a number, $b(x,r)$ a positive function
of $x \in M$, $r > 0$, and $h_\rho(x,y) := \rho$ or $d_M(x,y)$ for
$x,y \in M$, $\rho > 0$.  For a number $\rho > 0$ and a measurable map
$u : M \to Y$, we define the \emph{$\rho$-approximating energy density
  $e_u^\rho : M \to [\,0,+\infty\,)$} of $u$ by
\[
e_u^\rho(x) := \frac{1}{b(x,\rho)} \int_{B(x,\rho) \setminus \{x\}}
\left(\frac{d_Y(u(x),u(y))}{h_\rho(x,y)}\right)^p\;dy,
\qquad x \in M,
\]
and the \emph{$\rho$-approximating energy $E^\rho(u)$ of $u$} by
\[
E^\rho(u) := \frac{1}{2} \int_M e_u^\rho(x)\;dx
\in [\,0,+\infty\,],
\]
(see \cite{Jt:equilib, KvS:sobharm, St:heat, KwSy:sobmet}).
We assume the following two conditions:

\begin{itemize}
\item[(M)] There exists a positive function $\Theta(R)$, $R > 0$, with
  $\lim_{R \to 0} \Theta(R) = 1$ such that for any $\rho$, $R$ with $0
  < \rho \le R/2$, and for any measurable map $u : M \to Y$, we have
  \[
  e_u^R \le \Theta(R)\, e_u^\rho \quad\text{a.e. on $M$}.
  \]
\item[(AR)] There exists a constant $\kappa \in (\,0,1\,]$ such that
  \[
  \kappa\,b(x,r) \le |B(x,r)| \le b(x,r)
  \]
  for any $x \in M$ and $r > 0$.
\end{itemize}

In \cite{St:heat, KwSy:sobmet}, it is shown that some Bishop and
Bishop-Gromov type inequalities (called the measure contraction
property) imply the condition (M).  In particular, closed Riemannian
manifolds with volume measure and compact Alexandrov spaces with
Hausdorff measure all satisfy (M) and (AR) for a suitable function
$b(x,r)$ (see \cite{KwSy:sobmet}).  It follows from (M) that, as $\rho
\to 0$, $E^\rho$ $\Gamma$-converges to a lower semi-continuous
functional $E$ on the set of measurable maps from $M$ to $Y$ with
respect to $d_{L^p}$ (see \cite{St:heat, KwSy:sobmet}).  We call the
$\Gamma$-limit $E$ the \emph{energy functional}.  If $Y = \R$, then
$E^\rho$ and $E$ are both quadratic forms.  If $M$ and $Y$ are
Riemannian, then, for a suitably chosen $b(x,r)$, $E$ is the usual
energy functional defined by their Riemannian metrics upto a constant
multiple.  If $M$ is Riemannian and if $Y$ is a general metric space,
then $E$ is what Korevaar and Schoen studied in \cite{KvS:sobharm}.

\begin{thm} \label{thm:approxenergy}
  Assume that the $1$-local covering order of $M$ is at most
  $o(r^{-p})$ and that {\rm(M)} and {\rm(AR)} are both satisfied.
  Then, the $\rho$-approximating energy functional $E^\rho$ compactly
  converges to the energy functional $E$ as $\rho \to 0$ with respect
  to $d_{L^p}$.
\end{thm}

\begin{proof}
  Since $E^\rho$ $\Gamma$-converges to $E$, it suffices to
  prove that $\{E^\rho\}$ is asymptotically compact.
  For any $r$ and $\rho$ with $0 < \rho \le r$, we have
  \begin{align*}
    &\frac{1}{B(x,r)} \iint_{B(x,r) \times B(x,r)} d_Y(u(y),u(z))^p
    \; dzdy\\
    &\le \frac{1}{B(x,r)} \int_{B(x,r)}
    \int_{B(y,2r)} d_Y(u(y),u(z))^p \; dzdy\\
    &\le \frac{2\,b(x,r)\,(2r)^p}{|B(x,r)|}
    \int_{B(x,r)} e_u^{2r}(y) dy\\
    &\le 2^{p+1} \kappa^{-1} \Theta(2r)\,r^p
    \int_{B(x,r)} e_u^\rho(y) dy,
  \end{align*}
  which verifies the Poincar\'e inequality $(P)_{p,1,C,\rho,R}$ for
  some $C$ and $R$.  Applying Theorem \ref{thm:Pasympcpt} yields the
  asymptotic compactness of $\{E^\rho\}$.  This completes the proof.
\end{proof}

Theorem \ref{thm:approxenergy} and Proposition \ref{prop:Ecpt}
imply the following:

\begin{cor}
  $E$ is compact with respect to $d_{L^p}$.
\end{cor}

Theorem \ref{thm:approxenergy} and Proposition \ref{prop:min}(2)
imply:

\begin{cor}
  Let $\rho_i \to 0$ be a sequence of positive numbers and let a net of
  measurable maps $u_i : M \to Y$ be a $d_{L^p}$-bounded asymptotic
  minimizer of $\{E^{\rho_i}\}$ (see Definition \ref{defn:min}).
  Then, $\{u_i\}$ has a subsequence $L^p$-converging to a minimizer of
  $E$.
\end{cor}

The following proposition says that $E^\rho$ is non-compact typically,
though $\{E^\rho\}$ is asymptotically compact.

\begin{prop} \label{prop:noncpt}
  Assume that $h_\rho(x,y) := \rho$ and let $y_0 \in Y$ be a point.
  Then, $E^\rho$ is compact for a number $\rho > 0$ iff
  $L^p_{y_0}(M,Y)$ is proper.
\end{prop}

\begin{proof}
  If $L^p_{y_0}(M,Y)$ is proper, the compactness of $E^\rho$ is
  trivial.
  
  Let us prove the converse.
  Take any measurable map $u : M \to Y$.
  Since
  \[
  d_Y(u(x),u(y))^p \le 2^{p-1}\,(d_Y(u(x),y_0)^p + d_Y(u(y),y_0)^p),
  \]
  we have
  \begin{align} \label{eq:noncpt}
    E^\rho(u) &= \frac{1}{2} \iint_{0 < d_M(x,y) < \rho}
    \frac{1}{b(x,\rho)} \frac{d_Y(u(x),u(y))^p}{\rho^p}\;dxdy\\
    &\le \frac{2^{p-1}}{2\rho^p \inf_{x \in M} b(x,\rho)}
    \iint_{M \times M} (d_Y(u(x),y_0)^p + d_Y(u(y),y_0)^p)\;dxdy \notag\\
    &\le \frac{2^{p-1}\,|M|}{\rho^p \inf_{x \in M} b(x,\rho)}
    \,d_{L^p}(u,y_0)^p. \notag
  \end{align}
  It follows from (AR) that $\inf_{x \in M} b(x,\rho) > 0$.
  Assume that $E^\rho$ is compact.
  Take an $L^p$-bounded net $u_i \in L^p_{y_0}(M,Y)$.
  Then, by \eqref{eq:noncpt}, $E^\rho(u_i)$ is bounded.
  The compactness of $E^\rho$ tells us that $\{u_i\}$ has
  an $L^p$-convergent subnet.  This completes the proof.
\end{proof}

We now assume that $p = 2$ and $Y = \R$, so that $E^\rho$ and $E$ are
densely defined, nonnegative, and quadratic forms on $L^2(M)$.
Proposition \ref{prop:noncpt} and Theorems \ref{thm:approxenergy},
\ref{thm:cptspectra} imply:

\begin{cor}
  The spectrum for $E$ is discrete {\rm(}with only finite
  multiplicities{\rm)}.  Whenever $L^2(M)$ has infinite dimension,
  each $E^\rho$ is non-compact and has nonempty essential spectrum.
  As $\rho \to 0+$, the bottom of the essential spectrum for $E^\rho$
  is divergent to $+\infty$ and for any fixed $k \in \N$, the $k^{th}$
  eigenvalue for $E^\rho$ converges to that for $E$.  The
  eigenfunctions for $E^\rho$ also converge to eigenfunctions for $E$
  in the sense of Corollary \ref{cor:spectra}.
\end{cor}

\subsection{Compactness and convergence with a lower bound of Ricci curvature}

For constants $n \ge 2$ and $D > 0$, let $\mathcal{M}$ be the family
of $n$-dimensional closed Riemannian manifold of Ricci curvature $\ge
-(n-1)$ and diameter $\le D$.  We employ the probability measure
induced from the Riemannian volume measure on each $M \in
\mathcal{M}$.  The Gromov compactness theorem shows that the measured
Gromov-Hausdorff closure of $\mathcal{M}$ is compact (see (2.11) of
\cite{Fk:laplace}).  Let $\mathcal{Y}$ be a Gromov-Hausdorff compact
family of proper pointed metric spaces.  For $M \in \mathcal{M}$ and
$(Y,y) \in \mathcal{Y}$, we denote by $E_{M,Y} : L^2_y(M,Y) \to
[\,0,+\infty\,]$ the Korevaar-Schoen type energy (see
\cite{KvS:sobharm}).  Consider the family
$E_{\mathcal{M},\mathcal{Y}}$ of $E_{M,Y}$ for all $M \in \mathcal{M}$
and $Y \in \mathcal{Y}$.  We have the following:

\begin{thm} \label{thm:Ric}
  Any net of $E_{\mathcal{M},\mathcal{Y}}$ is asymptotically compact
  and has a compactly convergent subnet.
\end{thm}

\begin{proof}
  By a result of Buser \cite{Bs:isoper} we have a uniform bound of
  Poincar\'e constants, namely $(P)_{2,c,C,0,R}$ in Section
  \ref{ssec:Pcov} holds for $Y = \R$ and $E_{M,Y} \in
  E_{\mathcal{M},\mathcal{Y}}$.  Moreover, we have a uniform doubling
  constant by the Bishop-Gromov volume comparison theorem, which
  implies that any measured Gromov-Hausdorff limit of $\mathcal{M}$
  satisfies the doubling condition.  According to \cite{KST:Poincare},
  $(P)_{2,c,C,0,R}$ for $Y = \R$ implies $(P)_{2,c,C,0,R}$ for general
  $Y$.  Thus, the assumption of Theorem \ref{thm:Pasympcpt} is
  satisfied for any net of functionals in
  $E_{\mathcal{M},\mathcal{Y}}$.  This completes the proof.
\end{proof}

Consequently, by Theorem \ref{thm:cpt}, the family of $(\{u \in
L^2_y(M,Y) \mid E_{M,Y}(u) \le c\},y)$, $M \in \mathcal{M}$, $(Y,y)
\in \mathcal{Y}$, is relatively compact with respect to the pointed
Gromov-Hausdorff topology.  Let $\mathcal{H}$ be a Gromov-Hausdorff
compact family of proper pointed $\CAT(0)$-spaces.  Then, Proposition
\ref{prop:EJasympcpt} implies that for any net $\{E_i\}$ in
$E_{\mathcal{M},\mathcal{H}}$, the net of the resolvents $J^i_\lambda$
of $E_i$ is asymptotically compact.  By Theorem \ref{thm:cptconv},
$J^i_\lambda$ has a compactly convergent subnet.  Also, by
Propositions \ref{prop:Mosco}, \ref{prop:resolvent1}, and Theorem
\ref{thm:Gamma}, the Moreau-Yosida approximation $E_i^\lambda$ of
$E_i$ has a strongly convergent subnet.

Combining Theorem \ref{thm:Ric} and Kasue's result (Theorem 5.1 of
\cite{Ks:convLapII}) together with the work of Cheeger and Colding
\cite{CC:strRicIII} implies the following.

\begin{cor}
  Assume that $(Y,y)$ is a pointed complete Riemannian manifold.  Let
  $M_i \in \mathcal{M}$ be a net converging to a measured metric space
  with respect to the measured Gromov-Hausdorff topology.  Then,
  $E_{M_i,Y}$ compactly converges to the energy functional defined in
  \cite{Ks:convharm, Ks:convLapII}.
\end{cor}

Relative to this corollary, we refer Theorem 5.4 and Remark 5.1 of
\cite{KwSy:specstr} for the case $Y = \R$.

\begin{conj}
  Assume that a net $M_i \in \mathcal{M}$ measured Gromov-Hausdorff
  converges to a measured metric space $M$ and that a net $(Y_i,y_i)
  \in \mathcal{Y}$ pointed Gromov-Hausdorff converges to a proper
  metric space $(Y,y)$.
  The limit of $E_{M_i,Y_i}$ could uniquely determined only by
  $M$ and $(Y,y)$.
\end{conj}

\bibliographystyle{amsplain}
\bibliography{all}

\nocite{DM:Gamma}
\nocite{Gr:greenbook}
\nocite{Ms:compmedia}
\nocite{Jt:equilib, Jt:nlinDir, Jt:genharm}
\nocite{KwSy:specstr, KwSy:sobmet}
\nocite{Fk:laplace}

\end{document}